\documentclass[11pt,twoside]{article}
\newcommand{\ifims}[2]{#1} 
\newcommand{\ifAMS}[2]{#1}   
\newcommand{\ifau}[4]{#1}  
\newcommand{\ifbook}[2]{#1}   
\newcommand{\ifapp}[2]{#2}  
\newcommand{\iffourG}[2]{#1}  

\usepackage{tocloft}

\advance\cftsecnumwidth 0.3em\relax
\advance\cftsubsecnumwidth 0.3em\relax

\usepackage{ifthen}
\usepackage[ampersand]{easylist}
\usepackage{amsmath,amssymb,amsthm}
\usepackage[bb=stixtwo]{mathalpha}
\usepackage{mathtools}
\usepackage{natbib}
\usepackage{epsfig,graphicx}
\usepackage{comment}
\usepackage{color}
\usepackage{srcltx}
\usepackage[mathscr]{euscript}
\usepackage[math]{easyeqn}
\usepackage{etoolbox}
\usepackage{tocloft}
\usepackage{nccmath}
\usepackage{hyperref}
\hypersetup{
            colorlinks,
            linkcolor=hookersgreen,
            linktoc=blue,
            citecolor=ultramarine,
            urlcolor=black,
            filecolor=black
            }
\usepackage{nameref}
\usepackage{multirow}

\usepackage{accents}
\usepackage{mathrsfs}
\usepackage{bbm}
\usepackage{bm}
\usepackage{algorithm}
\usepackage{algpseudocode}

\ifims{
\textheight=22cm
\textwidth=14.8cm
\topmargin=0pt
\oddsidemargin=1.0cm
\evensidemargin=1.0cm
\linespread{1.3}

}{ 
}

\numberwithin{equation}{section}
\numberwithin{figure}{section}
\newcounter{example}[section]
\numberwithin{example}{section}
\newcounter{remark}[section]
\numberwithin{remark}{section}
\newtheorem{theorem}{Theorem}[section]

\newtheorem{lemma}[theorem]{Lemma}

\newtheorem{exmp}[example]{Example}
\newtheorem{rmrk}[remark]{Remark}
\newenvironment{example}{\begin{exmp}\rm}{\end{exmp}}
\newenvironment{remark}{\begin{rmrk}\rm}{\end{rmrk}}

\ifbook{
    \newcommand{\Chapter}[1]{\section{#1}}
    \newcommand{\Section}[1]{\subsection{#1}}
    \newcommand{\Subsection}[1]{\subsubsection{#1}}
    \def\Chname{Section }

  }
  {
    \newcommand{\Chapter}[1]{\chapter{#1}}
    \newcommand{\Section}[1]{\section{#1}}
    \newcommand{\Subsection}[1]{\subsection{#1}}
    \def\Chname{Chapter }

  }

\usepackage{mathrsfs}
\renewcommand{\(}{$\,}
\renewcommand{\)}{\,$}

\def\nquad{\hspace{-1cm}}
\def\mquad{\hspace{-0.5cm}}
\def\eqdef{\stackrel{\operatorname{def}}{=}}

\DeclareMathAlphabet{\mathbbmsl}{U}{bbm}{bx}{sl}

\DeclareMathSymbol{\Alpha}{\mathalpha}{operators}{"41}
\DeclareMathSymbol{\Beta}{\mathalpha}{operators}{"42}
\DeclareMathSymbol{\Epsilon}{\mathalpha}{operators}{"45}
\DeclareMathSymbol{\Zeta}{\mathalpha}{operators}{"5A}
\DeclareMathSymbol{\Eta}{\mathalpha}{operators}{"48}
\DeclareMathSymbol{\Iota}{\mathalpha}{operators}{"49}
\DeclareMathSymbol{\Kappa}{\mathalpha}{operators}{"4B}
\DeclareMathSymbol{\Mu}{\mathalpha}{operators}{"4D}
\DeclareMathSymbol{\Nu}{\mathalpha}{operators}{"4E}
\DeclareMathSymbol{\Omicron}{\mathalpha}{operators}{"4F}
\DeclareMathSymbol{\Rho}{\mathalpha}{operators}{"50}
\DeclareMathSymbol{\Tau}{\mathalpha}{operators}{"54}
\DeclareMathSymbol{\Chi}{\mathalpha}{operators}{"58}
\DeclareMathSymbol{\omicron}{\mathord}{letters}{"6F}

\newcommand{\cc}[1]{\mathscr{#1}}
\newcommand{\bb}[1]{\boldsymbol{#1}}

\DeclareFontFamily{U}{mathx}{\hyphenchar\font45}
\DeclareFontShape{U}{mathx}{m}{n}{
<5><6><7><8><9><10>
<10.95><12><14.4><17.28><20.74><24.88>
mathx10
}{}
\DeclareSymbolFont{mathx}{U}{mathx}{m}{n}
\DeclareFontSubstitution{U}{mathx}{m}{n}
\DeclareMathAccent{\widebar}{0}{mathx}{"73}

\renewcommand{\bar}[1]{\widebar{#1}}

\renewcommand{\tilde}[1]{\widetilde{#1}}

\makeatletter
\def\mathcenterto#1#2{\mathclap{\phantom{#1}\mathclap{#2}}\phantom{#1}}
\let\old@widetilde\widetilde
\def\widetildeto#1#2{\mathcenterto{#2}{\old@widetilde{\mathcenterto{#1}{#2\,}}}}
\let\old@widehat\widehat
\def\widehatto#1#2{\mathcenterto{#2}{\old@widehat{\mathcenterto{#1}{#2\,}}}}
\makeatother

\newcommand{\thankstitle}[1]{\ifthenelse{\equal{#1}{}}{}{\thanks{#1}}}
\newcommand{\thanksau}[1]{\ifthenelse{\equal{#1}{}}{}{\thanks{#1}}}

\newcommand{\aua}[6]
{\def\authora{#1}
\def\runauthora{#2}
\def\addressa{#3}
\def\emaila{#4}
\def\affiliationa{#5}
\def\thanksa{#6}}

\def\theauthors{
\ifau{ 
  \author{
    \authora
    \thanksau{\thanksa}
    \\[5.pt]
    \addressa \\
    \texttt{ \emaila}
  }
}
{  
  \author{
    \authora
    \thanksau{\thanksa}
    \\[5.pt]
    \addressa \\
    \texttt{ \emaila}
    \and
    \authorb
    \thanksau{\thanksb}
    \\[5.pt]
    \addressb \\
    \texttt{ \emailb}
  }
}
{   
  \author{
    \authora
    \thanksau{\thanksa}
    \\[5.pt]
    \addressa \\
    \texttt{ \emaila}
    \and
    \authorb
    \thanksau{\thanksb}
    \\[5.pt]
    \addressb \\
    \texttt{ \emailb}
    \and
    \authorc
    \thanksau{\thanksc}
    \\[5.pt]
    \addressc \\
    \texttt{ \emailc}
  }
} {   
  \author{
    \authora
    \thanksau{\thanksa}
    \\[5.pt]
    \addressa \\
    \texttt{ \emaila}
    \and
    \authorb
    \thanksau{\thanksb}
    \\[5.pt]
    \addressb \\
    \texttt{ \emailb}
    \and
    \authorc
    \thanksau{\thanksc}
    \\[5.pt]
    \addressc \\
    \texttt{ \emailc}
    \and
    \authord
    \thanksau{\thanksd}
    \\[5.pt]
    \addressd \\
    \texttt{ \emaild}
  }
}
}

\renewcommand{\Gamma}{\varGamma}
\renewcommand{\Pi}{\varPi}
\renewcommand{\Sigma}{\varSigma}
\renewcommand{\Delta}{\varDelta}
\renewcommand{\Lambda}{\varLambda}
\renewcommand{\Psi}{\varPsi}
\renewcommand{\Phi}{\varPhi}
\renewcommand{\Theta}{\varTheta}
\renewcommand{\Omega}{\varOmega}
\renewcommand{\Xi}{\varXi}
\renewcommand{\Upsilon}{\varUpsilon}

\def\arginf{\operatornamewithlimits{arginf}}

\def\argmin{\operatornamewithlimits{argmin}}

\def\av{\bm{a}}

\def\uv{\bm{u}}

\def\wv{\bm{w}}
\def\xv{\bm{x}}

\def\zv{\bm{z}}

\def\Av{\bm{A}}

\def\xiv{\bm{\xi}}

\newenvironment{mylist}{
  	\setlength{\parskip}{0pt}
  	\setlength{\parsep}{0pt}
\begin{itemize}
  	\setlength{\itemsep}{0pt}
    \setlength{\labelsep}{10pt}
    \setlength{\itemindent}{0em}}
{\end{itemize}
}

\usepackage{color}

\definecolor{blue(pigment)}{rgb}{0.2, 0.2, 0.6}
\definecolor{ultramarine}{rgb}{0.07, 0.04, 0.56}
\definecolor{darkspringgreen}{rgb}{0.09, 0.45, 0.27}
\definecolor{hookersgreen}{rgb}{0.0, 0.44, 0.0}
\definecolor{hgreen}{rgb}{0.0, 0.44, 0.0}
\definecolor{plum(traditional)}{rgb}{0.56, 0.27, 0.52}
\definecolor{purple(html/css)}{rgb}{0.5, 0.0, 0.5}
\definecolor{magenta(dye)}{rgb}{0.79, 0.08, 0.48}

\def\fh{h}
\def\fp{\fs^{\circ}}

\def\bvp{\bvn_{\circ}}
\def\upsvp{\upsv^{\circ}}

\def\hspm{\hspace{1pt}}

\def\AFN{\mathbb{Z}}
\def\BFN{\mathbb{B}}
\def\afv{\bb{s}}

\def\dmax{\kappa}

\def\bvn{\bb{\mu}}

\def\DFN{\DPN}

\def\AFN{\mathbb{U}}

\def\Tens{\mathcal{T}}

\def\Tensp{\Tens^{\circ}}

\def\hmax{\mathsf{c}}

\def\hL{h}

\def\dagg{\prime}

\def\Eta{\mathcal{H}}

\def\HL{\mathbbmsl{g}}

\def\dltw{\delta}

\def\dltwb{\omega}

\def\dltwu{\tau}
\def\dltwd{\dltw^{\dagg}}
\def\dltwbd{\dltwb^{\dagg}}

\def\II{\mathcal{I}}

\def\R{\mathbbmsl{R}}
\def\E{\mathbbmsl{E}}

\def\kappa{\varkappa}

\def\T{\top}

\def\Id{\mathbbmsl{I}}

\def\alp{\alpha}

\def\avn{\av}

\def\AFN{\mathbbmsl{U}}
\def\Avm{\bb{M}}
\def\AvmP{\bb{M}}
\def\AvmGP{\Avm_{\hspace{-1pt}\GP}}

\def\bias{\mathsf{b}}

\def\BBH{B}

\def\Bfv{\operatorname{B}}

\def\CA{\mathcal{A}}

\def\dimp{p}

\def\Eta{\cc{H}}

\def\fs{f}

\def\fn{g}

\def\fG{f_{\GP}}

\def\fGu{h}

\def\GP{G}

\def\IF{\mathbbmsl{F}}

\def\IFN{\IFL}

\def\IFP{\IFN_{\circ}}

\def\IFN{\IF}

\def\Kappa{\cc{K}}

\def\pent{\operatorname{pen}}

\def\rhot{t}
\def\rhos{\rho^{*}}

\def\rr{\mathtt{r}}

\def\rrd{\bar{\rr}}

\def\rrn{\rr}

\def\Tau{T}

\def\uvd{\uv^{\circ}}

\def\ups{\upsilon}
\def\upsv{\bb{\ups}}

\def\upsvd{\upsv^{\circ}}
\def\upsvs{\upsv^{*}}

\def\upsvn{\upsvd}
\def\upsn{\ups^{\circ}}

\def\ups{\upsilon}
\def\upsv{\bb{\ups}}

\def\UV{\mathcal{U}}

\def\UVz{\UV}

\def\Ups{\varUpsilon}

\def\wv{\bb{w}}

\def\xd{x^{\circ}}

    \renewcommand{\Section}[1]{\section{#1}}
    \renewcommand{\Subsection}[1]{\subsection{#1}}

\def\thetitle{Sharp bounds in perturbed smooth optimization}
\def\theruntitle{Perturbed convex optimization}

\def\theabstract{
This paper studies the problem of perturbed convex and smooth optimization.
The main results describe how the solution and the value of the problem change if the objective function is perturbed.
Examples include linear, quadratic, and smooth additive perturbations.
Such problems naturally arise in statistics and machine learning, stochastic optimization, stability and robustness analysis,
inverse problems, optimal control, etc.
The results provide accurate expansions for the difference between the solution of the original problem and its perturbed counterpart
with an explicit error term.
}

\def\kwdp{62F10,62E17}
\def\kwds{62J12}

\def\thekeywords{linear perturbation; smooth penalty; third-order and fourth-order expansions}

\def\thankstitle{}

\aua
{Vladimir Spokoiny}
{Spokoiny, V. }
{
    Weierstrass Institute Berlin,  
    HSE 
    and IITP RAS, Moscow,
    \\
    Mohrenstr. 39, 10117 Berlin, Germany
    }
{spokoiny@wias-berlin.de}
{Weierstrass-Institute and Humboldt University Berlin}
{
  Financial support by the German Research Foundation (DFG) through the Collaborative Research Center 388 ``Rough paths'' is gratefully acknowledged.
}

\begin{document}
\thispagestyle{empty}
{
\title{\thetitle}
\theauthors

\maketitle
\begin{abstract}
{\footnotesize \theabstract}
\end{abstract}

\ifAMS
    {\par\noindent\emph{AMS 2010 Subject Classification:} Primary \kwdp. Secondary \kwds}
    {\par\noindent\emph{JEL codes}: \kwdp}

} 

\tableofcontents

\Chapter{Introduction}
\label{Spertintr}

For a smooth convex function \( \fs(\cdot) \), consider an optimization problem
\begin{EQA}
	\upsvs
	&=&
	\argmin_{\upsv \in \Ups} \fs(\upsv) \, .
\label{8cicjwe6tdkcodlds}
\end{EQA}
Here \( \Ups \) is an open subset of \( \R^{\dimp} \) for \( \dimp < \infty \) and we focus on unconstraint optimization.
This paper studies the following question: how do \( \upsvs \) and \( \fs(\upsvs) \) change if the objective function 
\( \fs \) is perturbed in some special way?
We show how any smooth perturbation can be reduced to a linear one with 
\begin{EQA}[c]
	\fn(\upsv) = \fs(\upsv) + \langle \upsv, \Av \rangle .
\label{8cuc8c7ytw353wfgfff}
\end{EQA}
For \( \upsvn = \argmin_{\upsv} \fn(\upsv) \), we describe the differences \( \upsvn - \upsvs \) and 
\( \fn(\upsvn) - \fn(\upsvs) \) in terms of the Newton correction \( \IFN^{-1} \Av \) with \( \IFN = \nabla^{2} \fs(\upsvs) \).

\medskip
\par\textbf{Motivation 1: Statistical inference for SLS models.}
Most of statistical estimation procedures can be represented as ``empirical risk minimization''.
This includes least squares, maximum likelihood, minimum contrast, and many other methods.
Let \( L(\upsv) \) be  a random loss function (or empirical risk) and \( \E L(\upsv) \) its expectation, 
\( \upsv \in \Ups \subseteq \R^{\dimp} \), \( \dimp < \infty \).
Consider
\begin{EQA}[rcccl]
	\tilde{\upsv}
	&=&
	\argmin_{\upsv \in \Ups} L(\upsv)
	\, ;
	\qquad
	\upsvs
	&=&
	\argmin_{\upsv \in \Ups} \E \, L(\upsv) .
\end{EQA}
{Aim:}
describe the estimation loss \( \tilde{\upsv} - \upsvs \) and the prediction loss (excess)
\( L(\tilde{\upsv}) - L(\upsvs) \).
\cite{OsBa2021} provides some finite sample accuracy guarantees under self-concordance of \( L(\upsv) \).
\cite{Sp2024} introduced a special class of \emph{stochastically linear smooth} (SLS) models 
for which the stochastic term \( \zeta(\upsv) = L(\upsv) - \E L(\upsv) \) linearly depends on \( \upsv \)
and thus, its gradient \( \nabla \zeta \) does not depend on \( \upsv \):
\(
	\zeta(\upsv) - \zeta(\upsvd)
	=
	\langle \upsv - \upsvd, \nabla \zeta \rangle 
\).
This means that \( L(\upsv) \) is a linear perturbation of \( \fs(\upsv) = \E L(\upsv) \).
This reduces the study of the maximum likelihood estimator \( \tilde{\upsv} \) to linear perturbation analysis.

\medskip

\par\textbf{Motivation 2: Uncertainty quantification by smooth penalization.}
Given a smooth penalty \( \pent_{\GP}(\upsv) \) indexed by parameter \( \GP \), define 
\begin{EQA}
	\fG(\upsv) 
	&=& 
	\fs(\upsv) + \pent_{\GP}(\upsv) \, .
\label{LGLum2Gu2g}
\end{EQA}
A leading example is given by a quadratic penalty \( \| \GP \upsv \|^{2}/2 \): 
\begin{EQA}
	\fG(\upsv) 
	&=& 
	\fs(\upsv) + \| \GP \upsv \|^{2}/2 .
\label{LGLum2Gu22}
\end{EQA}
Again, the question under study is the corresponding change of the solution \( \upsvs_{\GP} = \argmin_{\upsv} \fG(\upsv) \) 
and the value \( \fG(\upsvs_{\GP}) \) of this problem. 
Adding a quadratic penalty changes the curvature of the objective function but does not change its smoothness.
It is important to establish an expansion whose remainder is only sensitive to the smoothness of \( \fs \).
Similar questions arise in 
Bayesian inference and high-dimensional Laplace approximation; see e.g. 
\cite{stuart2010inverse},
\cite{nickl2022bayesian},
\cite{katsevich2024laplaceapproximationaccuracyhigh}, and references therein.
For applications to Bayesian variational inference; see \cite{KaRi2022}.


\medskip
\par\textbf{Motivation 3: Quasi-Newton iterations.} \,
Consider an optimization routine for \eqref{8cicjwe6tdkcodlds} delivering \( \upsv_{k} \) as a solution after the step \( k \).
We want to measure the accuracy \( \upsv_{k} - \upsvs \) and the deficiency \( \fs(\upsv_{k}) - \fs(\upsvs) \).
Define 
\begin{EQA}[c]
	\fs_{k}(\upsv) \eqdef \fs(\upsv) - \langle \upsv - \upsv_{k}, \nabla \fs(\upsv_{k}) \rangle .
\label{cvivfedud6edfv6efuy}
\end{EQA}
Then \( \nabla \fs_{k}(\upsv_{k}) = 0 \).
As this function is strongly convex, it holds \( \upsv_{k} = \argmin_{\upsv} \fs_{k}(\upsv) \),
while \( \fs \) is a linear pertirbation of \( \fs_{k} \) with \( \Av = \nabla \fs(\upsv_{k}) \).
The obtained results describe \( \upsvs - \upsv_{k} \) in terms of \( \nabla \fs(\upsv_{k}) \) and \( \nabla^{2} \fs(\upsv_{k}) \).
Similar derivations for tensor methods in convex optimization and further references can be found in
\cite{doikov2021minimizing,doikov2022local,doikov2024gradient} and
\cite{rodomanov2022rates,rodomanov2021new}.
An accurate estimate of \( \upsv_{k} - \upsvs \) and of \( \fs_{k}(\upsv_{k}) - \fs_{k}(\upsvs) \)
can be used 
for analysis of many iterative procedures in optimization and stochastic approximation; see 
\cite{polyak1992acceleration}, 
\cite{NeNe1994},
\cite{boyd2004convex},
\cite{nesterov2018lectures},
\cite{polyak2007newton},
among many others.

There exists a vast literature on sensitivity analysis in
perturbed optimization focusing on the asymptotic setup
with a small parametric perturbation of a non-smooth objective function in Hilbert/Banach spaces under non-smooth constraints; see e.g. 
\cite{bonnans2000perturbation},
\cite{Bonnans1996}, 
\cite{BoSh1998} among many other.
For geometric properties of perturbed optimization with applications to machine learning see 
\cite{berthet2020} and references therein.
For our motivating examples, the assumption of a small perturbation is too restrictive.
Instead, we limit ourselves to a finite-dimensional optimization setup and to a smooth perturbation.
This enables us to derive accurate non-asymptotic closed-form results.

\def\biasl{b}
\par\textbf{This paper contribution.} \,
The main results describe the leading terms and the remainder of the expansion 
for the solution of the perturbed optimization problem for the cases of a linear, quadratic, and smooth perturbation.
The accuracy of the expansion strongly depends on the smoothness of the perturbed objective function
given in terms of a metric tensor \( \DFN \) with \( \DFN^{2} \preceq \dmax^{2} \IFN \) for \( \IFN \eqdef \nabla^{2} \fs(\upsvs) \).
Under second-order smoothness, we show
\begin{EQA}
	\| \IFN^{1/2} (\upsvn - \upsvs + \IFN^{-1} \Av) \|
    & \leq &
    \frac{2 \sqrt{\dltwb}}{1 - \dmax^{2} \dltwb} \| \DFN \, \IFN^{-1} \Av \| 
    \, , 
    \\
    \bigl| 2 \fn(\upsvn) - 2 \fn(\upsvs) + \| \IFN^{-1/2} \Av \|^{2} \bigr|
    & \leq &
    \frac{\dltwb}{1 - \dmax^{2} \dltwb} \| \DFN \, \IFN^{-1} \Av \|^{2} \, ,
\end{EQA}
where \( \| \cdot \| \) means the standard Euclidean norm 
and \( \dltwb \) describes the relative error of the second-order Taylor approximation;
see Theorem~\ref{PFiWigeneric}.
Further, a third-order self-concordance condition with a small constant \( \dltwu_{3} \) helps to substantially improve the accuracy:
with \( \biasl \eqdef \| \DFN \, \IFN^{-1} \Av \| \)
\begin{EQA}[rcl]
    \| \DFN^{-1} \IFN (\upsvn - \upsvs + \IFN^{-1} \Av) \|
    & \leq &
    \mfrac{8}{9} \, \dltwu_{3} \biasl^{2} 
    \, ,
    \\
    \bigl| 2 \fn(\upsvn) - 2 \fn(\upsvs) + \| \IFN^{-1/2} \Av \|^{2} \bigr|
    & \leq &
    \mfrac{1}{2} \, \dltwu_{3} \biasl^{3} \, ;
    \qquad
\end{EQA}
see Theorems~\ref{Pconcgeneric2} and \ref{PFiWigeneric2}.
All these bounds are non-asymptotic, the remainder is given in a closed form 
via the norm of the Newton correction \( \IFN^{-1} \Av \) and it is nearly optimal.
The results only involve the smoothness constant \( \dltwu_{3} \).
Even more, the fourth-order self-concordance yields for \( \biasl = \| \DFN \, \IFN^{-1} \Av \| \)
\begin{EQA}[rcl]
    \| \DFN^{-1} \IFN \{\upsvn - \upsvs + \IFN^{-1} \Av + \IFN^{-1} \nabla \Tens(\IFN^{-1} \Av) \} \|
    & \leq &
    \Bigl( \mfrac{1}{6} \dltwu_{4} + \mfrac{28}{27} \dmax^{2} \dltwu_{3}^{2} \Bigr) \, \biasl^{3} 
    \, ,
	\\
	\Bigl| \fn(\upsvn) - \fn(\upsvs) + \frac{1}{2} \| \IFN^{-1/2} \Av \|^{2} + \Tens(\IFN^{-1} \Av) \Bigr|
    & \leq &
    \Bigl( \mfrac{1}{24} \dltwu_{4} + \mfrac{52}{243} \dmax^{2} \dltwu_{3}^{2} \Bigr) \, \biasl^{4} 
    \, ,
    \qquad
\end{EQA}
with the third-order tensor \( \Tens = \mfrac{1}{6} \nabla^{3} \fs(\upsvs) \); see Theorem~\ref{Pconcgeneric4}.
The skewness correction term \( \IFN^{-1} \nabla \Tens(\IFN^{-1} \Av) \) naturally arises in variational inference and 
sharp Laplace approximation \cite{katsevich2024laplaceapproximationaccuracyhigh}.
We also explain how the case of any smooth perturbation can be reduced to a linear one. 

 
\par\textbf{Organization of the paper.} \,
Section~\ref{Ssmoothness} introduces some concepts of local smoothness in the sense of the self-concordance condition from
\cite{NeNe1994}.
Section~\ref{Squadnquad} collects the main results on linearly perturbed optimization. 
The obtained results are extended to the case of quadratic perturbation in Section~\ref{Slinquadr} 
and of a smooth perturbation in Section~\ref{Slinsmooth}.

\def\DFN{D}
\def\IFN{\mathbbmsl{F}}
\def\BFN{B}


%
\Section{Gateaux smoothness and self-concordance}
\label{Ssmoothness}
Below we assume \( \Ups \) to be an open convex domain in \( \R^{\dimp} \), \( \dimp < \infty \), and
the function \( \fs(\upsv) \), \( \upsv \in \Ups \subseteq \R^{\dimp} \) to be strongly convex with a positive definite Hessian 
\( \IFN(\upsv) \eqdef \nabla^{2} \fs(\upsv) \).
Also, assume \( \fs(\upsv) \) three or sometimes even four times Gateaux differentiable in \( \upsv \in \Ups \).
Local smoothness of \( \fs \) will be described by the relative error of the Taylor expansion 
of the third or fourth order.
Define
\begin{EQ}[rcl]
	\dltw_{3}(\upsv,\uv) 
	&=& 
	\fs(\upsv + \uv) - \fs(\upsv) - \langle \nabla \fs(\upsv), \uv \rangle 
	- \frac{1}{2} \langle \nabla^{2} \fs(\upsv), \uv^{\otimes 2} \rangle , 
	\\
	\dltwd_{3}(\upsv,\uv) 
	&=&
	\langle \nabla \fs(\upsv + \uv), \uv \rangle - \langle \nabla \fs(\upsv), \uv \rangle 
	- \langle \nabla^{2} \fs(\upsv), \uv^{\otimes 2} \rangle \, ,
\label{dltw3vufuv12f2ga}
\end{EQ}
and
\begin{EQA}
	\dltw_{4}(\upsv,\uv)
	& = &
	\fs(\upsv + \uv) - \fs(\upsv) - \langle \nabla \fs(\upsv), \uv \rangle 
	- \frac{1}{2} \langle \nabla^{2} \fs(\upsv), \uv^{\otimes 2} \rangle
	- \frac{1}{6} \langle \nabla^{3} \fs(\upsv), \uv^{\otimes 3} \rangle \, .
\label{hvcduywgedfuyg2y1y35e3wweg}
\end{EQA}
Now, for each \( \upsv \), suppose to be given a positive symmetric matrix 
\( \DFN(\upsv) \) 
defining a local metric and a local vicinity around \( \upsv \): for some radius \( \rr(\upsv) \),
\begin{EQA}
	\UVz(\upsv)
	&=&
	\bigl\{ \uv \in \R^{\dimp} \colon \| \DFN(\upsv) \uv \| \leq \rr(\upsv) \bigr\}
\label{ed7sycf7wedwgedq2ftwdfgtv}
\end{EQA}
Local smoothness properties of \( \fs \) at \( \upsv \) are given via the quantities
\begin{EQA}[rcccl]
    \dltwb(\upsv)
    & \eqdef &
    \sup_{\uv \colon \| \DFN(\upsv) \uv \| \leq \rr(\upsv)} \,
    \frac{2|\dltw_{3}(\upsv,\uv)|}{\| \DFN(\upsv) \uv \|^{2}} 
    \,\, ,
    \quad
    \dltwbd(\upsv)
    & \eqdef &
    \sup_{\uv \colon \| \DFN(\upsv) \uv \| \leq \rr(\upsv)} \, \frac{|\dltwd_{3}(\upsv,\uv)|}{\| \DFN(\upsv) \uv \|^{2}} \,\, . 
    \qquad
\label{dtb3u1DG2d3GPg}
\end{EQA}
The definition yields for any \( \uv \) with \( \| \DFN(\upsv) \uv \| \leq \rr \)
\begin{EQ}[rcccl]
	\bigl| \dltw_{3}(\upsv,\uv) \bigr|
	& \leq &
	\frac{\dltwb(\upsv)}{2} \| \DFN(\upsv) \uv \|^{2} 
	\, ,
	\qquad
	\bigl| \dltwd_{3}(\upsv,\uv) \bigr|
	& \leq &
	\dltwbd(\upsv) \| \DFN(\upsv) \uv \|^{2}
	\, .
	\qquad
\label{dta3u1DG2d3GPa1g}
\end{EQ}
%
The approximation results can be improved 
under a third-order upper bound on the error of Taylor expansion. 

\begin{description}
    \item[\label{LL3tref} \( \bb{(\mathcal{T}_{3})} \)]
      \textit{For some \( \dltwu_{3} \)}
\begin{EQA}
	\bigl| \dltw_{3}(\upsv,\uv) \bigr|
	& \leq &
	\frac{\dltwu_{3}}{6} \| \DFN(\upsv) \, \uv \|^{3} \, ,
	\quad
	\bigl| \dltwd_{3}(\upsv,\uv) \bigr|
	\leq 
	\frac{\dltwu_{3}}{2} \| \DFN(\upsv) \, \uv \|^{3} \, ,
	\quad
	\uv \in \UVz(\upsv).
\label{bd3xu16f3uo3st}
\end{EQA}
\end{description}
 
\begin{description}
    \item[\label{LL4tref} \( \bb{(\mathcal{T}_{4})} \)]
      \textit{For some \( \dltwu_{4} \)}
\begin{EQA}
	\bigl| \dltw_{4}(\upsv,\uv) \bigr|
	& \leq &
	\frac{\dltwu_{4}}{24} \| \DFN(\upsv) \, \uv \|^{4} \, ,
	\qquad
	\uv \in \UVz(\upsv).
\label{1mffmxum5st}
\end{EQA}
\end{description}

We also present a version of \nameref{LL3tref} resp. \nameref{LL4tref} in terms of the third (resp. fourth) derivative of \( \fs \).
\begin{description}
    \item[\label{LLsT3ref} \( \bb{(\mathcal{T}_{3}^{*})} \)]
    \emph{\( \fs(\upsv) \) is three times differentiable and 
	}
\begin{EQA}
    \sup_{\uv \colon \| \DFN(\upsv) \uv \| \leq \rr(\upsv)} \,\, \sup_{\zv \in \R^{\dimp}} \,\, 
    \frac{\bigl| \langle \nabla^{3} \fs(\upsv + \uv), \zv^{\otimes 3} \rangle \bigr|}
		 {\| \DFN(\upsv) \zv \|^{3}} 
	& \leq &
	\dltwu_{3} \, .
\label{jcxhydtferyu9j3d6vhew}
\end{EQA}

    \item[\label{LLsT4ref} \( \bb{(\mathcal{T}_{4}^{*})} \)]
    \emph{\( \fs(\upsv) \) is four times differentiable and 
	}
\begin{EQA}
    \sup_{\uv \colon \| \DFN(\upsv) \uv \| \leq \rr(\upsv)} \,\, \sup_{\zv \in \R^{\dimp}} \,\, 
    \frac{\bigl| \langle \nabla^{4} \fs(\upsv + \uv), \zv^{\otimes 4} \rangle \bigr|}
		 {\| \DFN(\upsv) \zv \|^{4}} 
	& \leq &
	\dltwu_{4} \, .
\label{jcxhydtferyu9j3d6vhew4}
\end{EQA}

\end{description}

%
\noindent
By Banach's characterization \cite{Banach1938}, \nameref{LLsT3ref} implies
\begin{EQA}
	\bigl| \langle \nabla^{3} \fs(\upsv + \uv), \zv_{1} \otimes \zv_{2} \otimes \zv_{3} \rangle \bigr|
	& \leq &	 
	\dltwu_{3} \| \DFN(\upsv) \zv_{1} \| \, \| \DFN(\upsv) \zv_{2} \| \, \| \DFN(\upsv) \zv_{3} \| \, 
\label{jbuyfg773jgion94euyyfg}
\end{EQA}
for any \( \uv \) with \( \| \DFN(\upsv) \uv \| \leq \rr \) and all \( \zv_{1} , \zv_{2}, \zv_{3} \in \R^{\dimp} \).
Similarly under \nameref{LLsT4ref}
\begin{EQA}
	&& \hspace{-.7cm}
	\bigl| \langle \nabla^{4} \fs(\upsv + \uv), \zv_{1} \otimes \zv_{2} \otimes \zv_{3} \otimes \zv_{4} \rangle \bigr|
	\\
	& \leq &
	\dltwu_{4} \| \DFN(\upsv) \zv_{1} \| \, \| \DFN(\upsv) \zv_{2} \| \, \| \DFN(\upsv) \zv_{3} \| \, \| \DFN(\upsv) \zv_{4} \| \, ,
	\quad 
	\zv_{1} , \zv_{2}, \zv_{3}, \zv_{4} \in \R^{\dimp} \, .
	\qquad
\label{jbuyfg773jgion94euyyfg4}
\end{EQA}

\begin{lemma}
\label{LdltwLa3t}
Under \nameref{LL3tref} or \nameref{LLsT3ref},
it holds for \( \dltwb(\upsv) \) and \( \dltwbd(\upsv) \) from \eqref{dtb3u1DG2d3GPg}
\begin{EQA}[rcccl]
\label{gtcdsftdffrvsewsea}
	\dltwb(\upsv)
	& \leq &
	\frac{\dltwu_{3} \, \rr(\upsv)}{3 } \, ,
	\qquad
	\dltwbd(\upsv)
	& \leq &
	\frac{\dltwu_{3} \, \rr(\upsv)}{2} 
	\, .
\label{gtcdsftdfvtwdsefhfdvfrvsewseG}
\end{EQA}
\end{lemma}

\begin{proof}
For any \( \uv \in \UVz(\upsv) \) with \( \| \DFN(\upsv) \uv \| \leq \rr \)
\begin{EQA}
	\bigl| \dltw_{3}(\upsv,\uv) \bigr|
	& \leq &
	\frac{\dltwu_{3}}{6} \, \| \DFN(\upsv) \uv \|^{3} 
	\leq 
	\frac{\dltwu_{3} \, \rr(\upsv)}{6} \, \| \DFN(\upsv) \uv \|^{2},
\label{jrgeteteer2234587654}
\end{EQA}
and the bound for \( \dltwb(\upsv) \) follows.
The proof for \( \dltwbd(\upsv) \) is similar.
Under \nameref{LLsT3ref}, apply the third order Taylor expansion to the 
univariate function \( \fs(\upsv + t \uv) \) of \( t \) and 
\nameref{LLsT3ref} with \( \zv \equiv \uv \).
\end{proof}

The values \( \dltwu_{3} \) and \( \dltwu_{4} \) are usually very small.
Some quantitative bounds are given later in this section
under the assumption that the function \( \fs(\upsv) \) can be written in the form \( \fs(\upsv) = n \hL(\upsv) \) 
for a fixed smooth function \( h(\upsv) \) with the Hessian \( \nabla^{2} \hL(\upsv) \). 
The factor \( n \) has meaning of the sample size%
\ifapp{; see \Chname \ref{ScritdimMLE} or \Chname \ref{SGBvM}.}{.}

\begin{description}
    \item[\label{LLtS3ref} \( \bb{(\mathcal{S}_{3}^{*})} \)]
      \emph{ \( \fs(\upsv) = n \hL(\upsv) \) for \( \hL(\upsv) \) three times differentiable and for a metric tensor \( \HL(\upsv) \)
\begin{EQA}
	\sup_{\uv \colon \| \HL(\upsv) \uv \| \leq \rr(\upsv)/\sqrt{n}} 
	\frac{\bigl| \langle \nabla^{3} \hL(\upsv + \uv), \uv^{\otimes 3} \rangle \bigr|}{\| \HL(\upsv) \uv \|^{3}}
	& \leq &
	\hmax_{3} \, .
\end{EQA}
}
    \item[\label{LLtS4ref} \( \bb{(\mathcal{S}_{4}^{*})} \)]
      \emph{ the function \( \hL(\cdot) \) satisfies \nameref{LLtS3ref} and  
\begin{EQA}
	\sup_{\uv \colon \| \HL(\upsv) \uv \| \leq \rr(\upsv)/\sqrt{n}}
	\frac{\bigl| \langle \nabla^{4} \hL(\upsv + \uv), \uv^{\otimes 4} \rangle \bigr|}{\| \HL(\upsv) \uv \|^{4}}
	& \leq &
	\hmax_{4} \, .
\end{EQA}
}
\end{description}

\noindent
\nameref{LLtS3ref} and \nameref{LLtS4ref}
are local versions of the so-called self-concordance condition; see \cite{NeNe1994} and \cite{OsBa2021}.
In fact, they require that each univariate function \( \hL(\upsv + t \uv) \) of \( t \in \R \)
is self-concordant with some universal constants \( \hmax_{3} \) and \( \hmax_{4} \).
Under \nameref{LLtS3ref} and \nameref{LLtS4ref}, with \( \DFN^{2}(\upsv) = n \, \HL^{2}(\upsv) \), the values 
\( \dltw_{3}(\upsv,\uv) \), \( \dltw_{4}(\upsv,\uv) \), and \( \dltwb(\upsv) \), \( \dltwbd(\upsv) \) can be well
bounded.

\begin{lemma}
\label{LdltwLaGP}
Suppose \nameref{LLtS3ref}.
Then 
\nameref{LL3tref} follows with \( \dltwu_{3} = \hmax_{3} n^{-1/2} \).
Moreover, for \( \dltwb(\upsv) \) and \( \dltwbd(\upsv) \) from \eqref{dtb3u1DG2d3GPg}, it holds
\begin{EQA}[rcccl]
	\dltwb(\upsv)
	& \leq &
	\frac{\hmax_{3} \, \rr(\upsv)}{3 n^{1/2}} \, ,
	\qquad
	\dltwbd(\upsv)
	& \leq &
	\frac{\hmax_{3} \, \rr(\upsv)}{2 n^{1/2}} \, .
\label{gtcdsftdfvtwdsefhfdvfrvsewseGP}
\end{EQA}
Also \nameref{LL4tref} follows from \nameref{LLtS4ref} with \( \dltwu_{4} = \hmax_{4} n^{-1} \).
\end{lemma}

\begin{proof}
For any \( \uv \in \UVz(\upsv) \), by third order Taylor expansion,
\begin{EQA}
	|\dltw_{3}(\upsv,\uv)|
	& \leq &
	\sup_{t \in [0,1]}
	\frac{1}{6} \bigl| \langle \nabla^{3} \fs(\upsv + t \uv), \uv^{\otimes 3} \rangle \bigr|
	=
	\frac{n}{6} \, 
	\sup_{t \in [0,1]}
	\bigl| \langle \nabla^{3} \hL(\upsv + t \uv), \uv^{\otimes 3} \rangle \bigr|
	\\
	& \leq &
	\frac{n \, \hmax_{3}}{6} \, \| \HL(\upsv) \uv \|^{3} 
	=
	\frac{n^{-1/2} \, \hmax_{3}}{6} \, \| \DFN(\upsv) \uv \|^{3}
	\leq 
	\frac{n^{-1/2} \, \hmax_{3} \, \rr(\upsv)}{6} \, \| \DFN(\upsv) \uv \|^{2} \, .
\label{jrgeteteer2234587654}
\end{EQA}
This implies \nameref{LL3tref} as well as \eqref{gtcdsftdfvtwdsefhfdvfrvsewseGP}; see \eqref{dta3u1DG2d3GPa1g}.
The statement about \nameref{LL4tref} is similar.
\end{proof}

\def\BFN{B}
\def\Bfv{B}
\def\rell{\ell}
\def\rellG{\rell_{\GP}}
\def\rbias{\rho}
\def\biasl{b}
\def\biasls{\biasl_{*}}

\def\AFN{U}
\def\CAD{\mathcal{A}}
\def\Reminder{\mathcal{E}}
\def\upsvf{\upsv^{\dag}}
\def\uvn{\uv^{\circ}}

\Section{Linearly perturbed optimization}
\label{Squadnquad}

Let \( \fs(\upsv) \) be a smooth convex function, 
\begin{EQA}
	\upsvs
	&=&
	\arginf_{\upsv} \fs(\upsv),
	\qquad
	\fs(\upsvs)
	=
	\inf_{\upsv} \fs(\upsv) ,
	\qquad
	\IFN = \nabla^{2} \fs(\upsvs) \, .
\label{fg5hg3gf98tkj3dciryt}
\end{EQA}
We study how the minimizer and the minimum value change under a linear perturbation.
More precisely, let a function \( \fn(\upsv) \) satisfy for some vector \( \Av \)
\begin{EQA}
	\fn(\upsv) - \fn(\upsvs) 
	&=&
	\bigl\langle \upsv - \upsvs, \Av \bigr\rangle + \fs(\upsv) - \fs(\upsvs) .
\label{4hbh8njoelvt6jwgf09}
\end{EQA}
Define
\begin{EQA}
	\upsvn
	& \eqdef &
	\arginf_{\upsv} \fn(\upsv),
	\qquad
	\fn(\upsvn)
	=
	\inf_{\upsv} \fn(\upsv) .
\label{6yc63yhudf7fdy6edgehy} 
\end{EQA}
We derive bounds and expansions for the quantities \( \upsvn - \upsvs \) and
\( \fn(\upsvn) - \fn(\upsvs) \).
First, we consider the case of a quadratic function \( \fs \).

\begin{lemma}
\label{Pquadquad}
Let \( \fs(\upsv) \) be quadratic with \( \nabla^{2} \fs(\upsv) \equiv \IFN \) and let 
\( \fn(\upsv) \) be from \eqref{4hbh8njoelvt6jwgf09}.
Then  
\begin{EQA}
	\upsvn - \upsvs
	&=&
	- \IFN^{-1} \Av,
	\qquad
	\fn(\upsvn) - \fn(\upsvs)
	=
	- \frac{1}{2} \| \IFN^{-1/2} \Av \|^{2} .
\label{kjcjhchdgehydgtdtte35}
\end{EQA}
\end{lemma}

\begin{proof}
If \( \fs(\upsv) \) is quadratic, then under \eqref{4hbh8njoelvt6jwgf09}, \( \fn(\upsv) \) is quadratic as well
with \( \nabla^{2} \fn(\upsv) \equiv \IFN \).
This implies
\begin{EQA}
	\nabla \fn(\upsvs) - \nabla \fn(\upsvn)
	&=&
	\IFN (\upsvs - \upsvn) .
\label{dcudydye67e6dy3wujhds7}
\end{EQA}
Further, \eqref{4hbh8njoelvt6jwgf09} and \( \nabla \fs(\upsvs) = 0 \) yield \( \nabla \fn(\upsvs) = \Av \).
Together with \( \nabla \fn(\upsvn) = 0 \), this implies
\( \upsvn - \upsvs = - \IFN^{-1} \Av \).
The Taylor expansion of \( \fn \) at \( \upsvn \) yields by \( \nabla \fn(\upsvn) = 0 \)
\begin{EQA}
	\fn(\upsvs) - \fn(\upsvn)
	&=&
	\frac{1}{2} \| \IFN^{1/2} (\upsvn - \upsvs) \|^{2}
	=
	\frac{1}{2} \| \IFN^{-1/2} \Av \|^{2} 
\label{8chuctc44wckvcuedje}
\end{EQA}
and the assertion follows.
\end{proof}

\Subsection{Localization of the perturbed minimizer}
Let \( \fs \) satisfy \eqref{dtb3u1DG2d3GPg} at \( \upsvs \) with 
\( \DFN(\upsvs) = \DFN \) satisfying \( \DFN^{2} \leq \dmax^{2} \, \IFN \) for some \( \dmax > 0 \).
Equivalently, \( \dmax^{2} \IFN - \DFN^{2} \) is positive semidefinite.
The next result describes the concentration properties of \( \upsvn \) from \eqref{6yc63yhudf7fdy6edgehy} in a local elliptic set
\begin{EQA}[c]
	\{ \upsv \colon \| \IFN^{1/2} (\upsv - \upsvs) \| \leq \rr \} ,
\label{0cudc7e3jfuyvct6eyhgwe}
\end{EQA}
where \( \rr \) is slightly larger than \( \| \IFN^{-1/2} \Av \| \).

\begin{theorem}
\label{Pconcgeneric}
Let 
\begin{mylist}
\item
\( \fs(\upsv) \) be a strongly convex function with \( \fs(\upsvs) = \inf_{\upsv} \fs(\upsv) \)  
and \( \IFN = \nabla^{2} \fs(\upsvs) \);
\item
\( \fn(\upsv) \) and \( \fs(\upsv) \) be related by \eqref{4hbh8njoelvt6jwgf09} with some vector \( \Av \);
\item
the function \( \fs(\upsv) \) satisfy \eqref{dtb3u1DG2d3GPg} at \( \upsv = \upsvs \) 
with \( \DFN = \DFN(\upsvs) \) and \( \dltwbd \geq 0 \);
\item
for some \( \dmax > 0 \) 
\begin{EQA}
	\DFN^{2} \leq \dmax^{2} \, \IFN
	\, ,
	\qquad
	\dltwbd \dmax^{2} < 1
	\, ,
	\qquad
	\rr(\upsvs)
	& \geq & 
	\frac{\dmax}{1 - \dltwbd \dmax^{2}} \| \IFN^{-1/2} \Av \|
	\, .
	\qquad
\label{rrm23r0ut3ua}
\end{EQA}
\end{mylist}
Then the minimizer \( \upsvn \) defined in \eqref{6yc63yhudf7fdy6edgehy} satisfies
\begin{EQA}
	\| \IFN^{1/2} (\upsvn - \upsvs) \|  
	& \leq &
	\mfrac{1}{1 - \dltwbd \dmax^{2}} \| \IFN^{-1/2} \Av \| 
	\, , 
	\\
	\| \DFN (\upsvn - \upsvs) \|
	& \leq &
	\mfrac{\dmax}{1 - \dltwbd \dmax^{2}} \, \| \IFN^{-1/2} \Av \| 
	\, . 
\label{rhDGtuGmusGU0a}
\end{EQA}
\end{theorem}

\begin{proof}
Let \( \rrn \) be any radius satisfying
\begin{EQA}[c]
	\mfrac{1}{1 - \dltwbd \dmax^{2}} \| \IFN^{-1/2} \Av \|
	\leq 
	\rrn
	\leq 
	\mfrac{\rr(\upsvs)}{\dmax} 
	\, ,
\label{juiuikiu8f6ce33tofg0}
\end{EQA} 
and let \( \upsv \) be a point on the boundary of the local set from \eqref{0cudc7e3jfuyvct6eyhgwe}
satisfying \( \| \IFN^{1/2} (\upsv - \upsvs) \| \leq \rrn \).
We also write \( \uv = \upsv - \upsvs \).
The idea is to show that the derivative  \( \frac{d}{dt} \fn(\upsvs + t \uv) \) 
is positive for \( t > 1 \).
Then all stationary points of \( \fn(\upsv) \) satisfy \( \| \IFN^{1/2} (\upsv - \upsvs) \| \leq \rrn \).
We use the decomposition
\begin{EQA}
	\fn(\upsvs + \rhot \uv) - \fn(\upsvs)
	&=&
	\langle \Av, \uv \rangle \, \rhot 
	+ \fs(\upsvs + \rhot \uv) - \fs(\upsvs) 
	\, .
\label{LGtsGtuLGtsa}
\end{EQA}
With \( \fGu(t) = \fs(\upsvs + \rhot \uv) - \fs(\upsvs) + \langle \Av, \uv \rangle \, \rhot \), it holds
\begin{EQA}
	\frac{d}{d \rhot} \fs(\upsvs + \rhot \uv)
	&=&
	- \langle \Av, \uv \rangle + \fGu'(\rhot) .
\label{frddtLtGstua}
\end{EQA}
By definition of \( \upsvs \), it also holds \( \fGu'(0) = \langle \Av, \uv \rangle \).
The identity \( \nabla^{2} \fs(\upsvs) = \IFN \) yields \( \fGu''(0) = \| \IFN^{1/2} \uv \|^{2} \).
Bound \eqref{dta3u1DG2d3GPa1g} implies for \( | \rhot | \leq 1 \)
by \( \| \DFN \uv \| \leq \dmax \| \IFN^{1/2} \uv \| \)
\begin{EQA}
	\bigl| \fGu'(\rhot) - \fGu'(0) - \rhot \fGu''(0) \bigr|
	& \leq &
	|\rhot| \, \| \IFN^{1/2} \uv \|^{2} \, \dltwbd \, \dmax^{2}
	\, .
\label{fptfp0fpttfpp13a}
\end{EQA}
For \( \rhot = 1 \), we obtain by \eqref{rrm23r0ut3ua} and 
\( | \langle \Av, \uv \rangle \bigr| = | \langle \IFN^{-1/2} \Av, \IFN^{1/2} \uv \rangle \bigr|
	\leq \| \IFN^{-1/2} \Av \| \; \| \IFN^{1/2} \uv \| \)
\begin{EQA}
	\fGu'(1) 
	& \geq &
	\langle \Av, \uv \rangle + \| \IFN^{1/2} \uv \|^{2} - \| \IFN^{1/2} \uv \|^{2} \, \dltwbd \, \dmax^{2}
	\\
	& \geq &
	\| \IFN^{1/2} \uv \| \bigl\{ \rrn (1 -  \dltwbd \, \dmax^{2}) - \| \IFN^{-1/2} \Av \| \bigr\}
	\geq 
	0 
	\, .
\label{fp1fpp13d3rGa}
\end{EQA}
Moreover, strong convexity of \( \fGu(\rhot) \) implies that \( \fGu'(\rhot) - \fGu'(0) \) increases in 
\( \rhot \) for \( \rhot > 1 \).
%
This implies 
\begin{EQA}
	\frac{d}{dt} \fn(\upsvs + \rhot \uv) 
	& > &
	0 
	\, ,
	\qquad
	t > 1 
	\, ,
\label{ddtLGtstu33a}
\end{EQA}
and the assertion follows.
\end{proof}

\Subsection{{Second-order expansions}}
Theorem~\ref{Pconcgeneric} localizes the perturbed minimizer \( \upsvn = \arginf_{\upsv} \fn(\upsv) \)
in the local neighborhood \eqref{0cudc7e3jfuyvct6eyhgwe} of \( \upsvs \).
Using smoothness properties of \( \fn \) or, equivalently, of \( \fs \), in this vicinity gives 
explicit expansions for \( \upsvn - \upsvs \) and for \( \fn(\upsvn) - \fn(\upsvs) \).

\begin{theorem}
\label{PFiWigeneric}
Assume the conditions of Theorem~\ref{Pconcgeneric}.
Then with \( \dltwb \leq \dltwbd \) from \eqref{dtb3u1DG2d3GPg}
\begin{EQ}[rcl]
    - \frac{\dltwb}{1 - \dmax^{2} \hspm \dltwb} \| \DFN \, \IFN^{-1} \Av \|^{2}
    & \leq &
    2 \fn(\upsvn) - 2 \fn(\upsvs) 
    + \| \IFN^{-1/2} \Av \|^{2}
    \\
    & \leq &
    \dltwb \hspm \| \DFN \, \IFN^{-1} \Av \|^{2} \, .
    \qquad 
\label{3d3Af12DGttGa}
\end{EQ}
Also
\begin{EQ}[rcl]
    \| \IFN^{1/2} (\upsvn - \upsvs + \IFN^{-1} \Av) \|
    & \leq &
    \frac{2 \sqrt{\dltwb}}{1 - \dmax^{2} \hspm \dltwb} \, \| \DFN \, \IFN^{-1} \Av \| \, ,
\label{DGttGtsGDGm13rGa}
\end{EQ}
\end{theorem}

\begin{proof}
Let \( \rrn \) satisfy \eqref{juiuikiu8f6ce33tofg0} and \( \| \IFN^{1/2} (\upsv - \upsvs) \| \leq \rrn \).
By \eqref{dtb3u1DG2d3GPg}, 
\begin{EQA}
	\Bigl| 
		\fs(\upsv) - \fs(\upsvs) - \mfrac{1}{2} \| \IFN^{1/2} (\upsv - \upsvs) \|^{2} 
	\Bigr|
	& \leq &
	\mfrac{\dltwb}{2} \| \DFN (\upsv - \upsvs) \|^{2} .
\label{d3GrGELGtsG12}
\end{EQA}
Further, 
\begin{EQA}[rcl]
	&& \nquad
	\fn(\upsv) - \fn(\upsvs) + \mfrac{1}{2} \| \IFN^{-1/2} \Av \|^{2}
	\\
	&=&
	\bigl\langle \upsv - \upsvs, \Av \bigr\rangle
	+ \fs(\upsv) - \fs(\upsvs) + \mfrac{1}{2} \| \IFN^{-1/2} \Av \|^{2} 
	\\
	&=&
	\mfrac{1}{2} \bigl\| \IFN^{1/2} (\upsv - \upsvs) + \IFN^{-1/2} \Av \bigr\|^{2}
	+ \fs(\upsv) - \fs(\upsvs) - \mfrac{1}{2} \| \IFN^{1/2} (\upsv - \upsvs) \|^{2} .
	\qquad 
\label{12ELGuELusG}
\end{EQA}
As \( \| \IFN^{1/2} (\upsvn - \upsvs) \| \leq \rrn \) and \( \upsvn \) minimizes \( \fn(\upsv) \), 
we derive by \eqref{d3GrGELGtsG12}
\begin{EQA}
	&& \nquad
	\fn(\upsvn) - \fn(\upsvs) + \mfrac{1}{2} \| \IFN^{-1/2} \Av \|^{2}
	=
	\inf_{\upsv \colon \| \IFN^{1/2} (\upsvn - \upsvs) \| \leq \rrn} 
	\Bigl\{ 
		\fn(\upsv) - \fn(\upsvs) + \mfrac{1}{2} \| \IFN^{-1/2} \Av \|^{2} 
	\Bigr\}
	\\
	& \geq &
	\inf_{\upsv \colon \| \IFN^{1/2} (\upsvn - \upsvs) \| \leq \rrn} 
	\Bigl\{ 
		\mfrac{1}{2} \bigl\| \IFN^{1/2} (\upsv - \upsvs) + \IFN^{-1/2} \Av \bigr\|^{2} 
		- \mfrac{\dltwb}{2} \| \DFN (\upsv - \upsvs) \|^{2}
	\Bigr\} .
\label{d3G1212222B} 
\end{EQA}
Denote \( \uv = \IFN^{1/2} (\upsv - \upsvs) \), \( \xiv = \IFN^{-1/2} \Av \), and 
\( \BFN = \IFN^{-1/2} \, \DFN^{2} \, \IFN^{-1/2} \).
As \( \DFN^{2} \leq \dmax^{2} \IFN \), 
it holds \( \| \BFN \| \leq \dmax^{2} \) 
for the operator norm of \( \BFN \) and 
\begin{EQA}
	&& \mquad
	\inf_{\upsv \colon \| \IFN^{1/2} (\upsvn - \upsvs) \| \leq \rrn} 
	\bigl\{ \bigl\| \IFN^{1/2} (\upsv - \upsvs) + \IFN^{-1/2} \Av \bigr\|^{2} 
		- \dltwb \| \DFN (\upsv - \upsvs) \|^{2} \bigr\}
	\\
	&=&
	\inf_{\| \uv \| \leq \rr} \bigl\{ \| \uv + \xiv \|^{2} - \dltwb \, \uv^{\T} \BFN \uv \bigr\}
	\geq 
	- \xiv^{\T} \bigl\{ (\Id_{\dimp} - \dltwb \, \BFN)^{-1} - \Id_{\dimp} \bigr\} \xiv
	\geq 
	- \mfrac{\dltwb}{1 - \dmax^{2} \dltwb} \xiv^{\T} \BFN \, \xiv
\label{d7eneyf6g53geygywn}
\end{EQA}
with \( \Id_{\dimp} \) being the identity matrix in \( \R^{\dimp} \).
This yields by \( \xiv^{\T} \BFN \xiv = \| \DFN \, \IFN^{-1} \Av \|^{2} \)
\begin{EQA}
	\fn(\upsvn) - \fn(\upsvs) + \mfrac{1}{2} \| \IFN^{-1/2} \Av \|^{2}
	& \geq &
	- \mfrac{\dltwb}{2(1 - \dmax^{2} \dltwb)} \| \DFN \, \IFN^{-1} \Av \|^{2} . 
\label{fd3G122B2} 
\end{EQA}
Similarly, using \( \| (\Id_{\dimp} + \dltwb \, \BFN)^{-1} \xiv \| \leq \| \xiv \| \leq \rrn \), we derive
\begin{EQA}
	&& \nquad 
	\fn(\upsvn) - \fn(\upsvs) + \mfrac{1}{2} \| \IFN^{-1/2} \Av \|^{2}
	\\
	& \leq &
	\inf_{\upsv \colon \| \IFN^{1/2} (\upsvn - \upsvs) \| \leq \rrn} 
	\Bigl\{ 
		\mfrac{1}{2} \bigl\| \IFN^{1/2} (\upsv - \upsvs) + \IFN^{-1/2} \Av \bigr\|^{2} 
		+ \mfrac{\dltwb}{2} \| \DFN (\upsv - \upsvs) \|^{2}
	\Bigr\}
	\\
	& \leq &
	\mfrac{1}{2} \xiv^{\T} \bigl\{ \Id_{\dimp} - (\Id_{\dimp} + \dltwb \, \BFN)^{-1} \bigr\} \xiv
	\leq 
	\mfrac{\dltwb }{2} \, \xiv^{\T} \BFN \, \xiv
	=
	\mfrac{\dltwb }{2} \, \| \DFN \, \IFN^{-1} \Av \|^{2} 
	\, . 
	\quad 
\label{fd3G122B2m} 
\end{EQA}
These bounds imply 
\eqref{3d3Af12DGttGa}.

Now we derive similarly to \eqref{12ELGuELusG} that for any \( \upsv \)
satisfying \( \| \IFN^{1/2} (\upsv - \upsvs) \| \leq \rrn \)
\begin{EQA}
	\fn(\upsv) - \fn(\upsvs) 
	& \geq &
	\bigl\langle \upsv - \upsvs, \Av \bigr\rangle
	+ \mfrac{1}{2} \| \IFN^{1/2} (\upsv - \upsvs) \|^{2}
	- \mfrac{\dltwb}{2} \| \DFN (\upsv - \upsvs) \|^{2} .
\label{LGvLGvsGf1d3G2}
\end{EQA}
A particular choice \( \upsv = \upsvn \) yields
\begin{EQA}
	\fn(\upsvn) - \fn(\upsvs) 
	& \geq &
	\bigl\langle \upsvn - \upsvs, \Av \bigr\rangle
	+ \mfrac{1}{2} \| \IFN^{1/2} (\upsvn - \upsvs) \|^{2}
	- \mfrac{\dltwb}{2} \| \DFN (\upsvn - \upsvs) \|^{2} .
\label{21GsvtvGDG3G2}
\end{EQA}
Combining this inequality with \eqref{fd3G122B2m} gives
\begin{EQA}
	\bigl\langle \upsvn - \upsvs, \Av \bigr\rangle
	+ \mfrac{1}{2} \| \IFN^{1/2} (\upsvn - \upsvs) \|^{2}
	- \mfrac{\dltwb}{2} \| \DFN (\upsvn - \upsvs) \|^{2} 
	& \leq &
	- \mfrac{1}{2} \xiv^{\T} (\Id_{\dimp} - \dltwb \, \BFN) \xiv .
\label{2m1DGd3G123G}
\end{EQA}
With 
\( \uvd = \IFN^{1/2} (\upsvn - \upsvs) \), this implies
\begin{EQA}
	2 \bigl\langle \uvd, \xiv \bigr\rangle + {\uvd}^{\T}(\Id_{\dimp} - \dltwb \BFN) \uvd 
	& \leq &
	- \xiv^{\T} (\Id_{\dimp} - \dltwb \, \BFN) \xiv \, ,
\label{dtxi2fd1d22}
\end{EQA}
and hence, by \( \| \BFN \| \leq \dmax^{2} \) and \( (\Id_{\dimp} - \dltwb \, \BFN)^{-1} \leq (1 - \dmax^{2} \dltwb)^{-1} \)
\begin{EQA}
	&& \nquad
	\bigl\{ \uvd + (\Id_{\dimp} - \dltwb \BFN)^{-1} \xiv \bigr\}^{\T} (\Id_{\dimp} - \dltwb \BFN) \bigl\{ \uvd + (\Id_{\dimp} - \dltwb \BFN)^{-1} \xiv \bigr\}
	\\
	& \leq &
	\xiv^{\T} \bigl\{ (\Id_{\dimp} - \dltwb \, \BFN)^{-1} - (\Id_{\dimp} - \dltwb \, \BFN) \bigr\} \xiv
	=
	\dltwb \xiv^{\T} \BFN \bigl\{ (\Id_{\dimp} - \dltwb \, \BFN)^{-1} + \Id_{\dimp} \bigr\} \xiv
	\\
	& \leq &
	\frac{2 \dltwb - \dmax^{2} \, \dltwb^{2}}{1 - \dmax^{2} \dltwb} \, \xiv^{\T} \BFN \xiv 
	\, ,
\label{uv11wxi22w1w}
\end{EQA}
yielding
\begin{EQA}[c]
	\| \uvd + (\Id_{\dimp} - \dltwb \BFN)^{-1} \xiv \|^{2}
	\leq 
	\frac{2 \dltwb - \dmax^{2} \, \dltwb^{2}}{(1 - \dmax^{2} \dltwb)^{2}} \, \xiv^{\T} \BFN \, \xiv
	\, .
\end{EQA}
Further, by \( \BFN \leq \dmax^{2} \Id_{\dimp} \),
\begin{EQA}[c]
	\| \xiv - (\Id_{\dimp} - \dltwb \BFN)^{-1} \xiv \|
	=
	\dltwb \| (\Id_{\dimp} - \dltwb \BFN)^{-1} \BFN \xiv \|
	\leq 
	\frac{\dltwb}{1 - \dmax^{2} \dltwb} \| \BFN \xiv \|
	\leq 
	\frac{\dmax \, \dltwb}{1 - \dmax^{2} \dltwb} \| \BFN^{1/2} \xiv \|
	\, .
\end{EQA}
The use of the triangle inequality, \( \dmax^{2} \dltwb < 1 \), and the inequality 
\( \alp + \sqrt{2 \alp - \alp^{2}} \leq 2 \sqrt{\alp} \) for \( \alp = \dmax^{2} \dltwb \) yields
\begin{EQA}[c]
	\| \uvd + \xiv \|
	\leq 
	\frac{\dmax \, \dltwb + \sqrt{2 \dltwb - \dmax^{2} \dltwb^{2}}}{1 - \dmax^{2} \dltwb} \, \| \BFN^{1/2} \xiv \|
	\leq 
	\frac{2 \sqrt{\dltwb}}{1 - \dmax^{2} \dltwb} \, \| \DFN \IFN^{-1} \Av \|
	\, ,
\end{EQA}
as stated.
\end{proof}


\Subsection{{Expansions under third-order smoothness}}
The results of Theorem~\ref{PFiWigeneric} can be refined if
\( \fs \) satisfies condition \nameref{LL3tref}.
Introduce the local set in the \( \DFN \)-metric
\begin{EQA}[rcl]
	\CAD(\rr)
	& \eqdef &
	\{ \upsv \colon \| \DFN (\upsv - \upsvs) \| \leq \rr \} 
	\, .
\label{0cudc7e3jfuyvct6eyhgweD}
\end{EQA}

\begin{theorem}
\label{Pconcgeneric2}
Let 
\begin{mylist}
\item
\( \fs(\upsv) \) be a strongly convex function with \( \fs(\upsvs) = \inf_{\upsv} \fs(\upsv) \)  
and \( \IFN = \nabla^{2} \fs(\upsvs) \);
\item
\( \fn(\upsv) \) and \( \fs(\upsv) \) be related by \eqref{4hbh8njoelvt6jwgf09} with some vector \( \Av \);
\item
the function \( \fs(\upsv) \) satisfy \nameref{LL3tref} at \( \upsvs \) with 
\( \DFN^{2} \), \( \rrn = \rr(\upsvs) \), and \( \dltwu_{3} \); 

\item
for some \( \dmax > 0 \) 
\begin{EQA}
	\DFN^{2} \leq \dmax^{2} \hspm \IFN , 
	\qquad 
	\rrn \geq \mfrac{4\dmax}{3} \, \| \IFN^{-1/2} \Av \| ,
	\qquad 
	\dmax^{3} \dltwu_{3} \, \| \IFN^{-1/2} \Av \|
	& \leq &
	\mfrac{3}{8} 
	\, .
	\qquad
\label{yxdhewndu7jwnjjuu}
\end{EQA}
\end{mylist}
Then \( \upsvn = \arginf_{\upsv} \fn(\upsv) \) satisfies with \( \rhos \eqdef \dmax^{3} \dltwu_{3} \, \| \IFN^{-1/2} \Av \| \)
\begin{EQA}
	\| \IFN^{1/2} (\upsvn - \upsvs) \|  
	& \leq &
	\mfrac{1}{1 - 2 \rhos/3} \| \IFN^{-1/2} \Av \| 
	\leq 
	\mfrac{4}{3} \| \IFN^{-1/2} \Av \|
	\, ,
	\\
	\| \DFN (\upsvn - \upsvs) \|  
	& \leq &
	\mfrac{\dmax}{1 - 2 \rhos/3} \, \| \IFN^{-1/2} \Av \| 
	\leq 
	\mfrac{4\dmax}{3} \| \IFN^{-1/2} \Av \|
	\, . 
\label{rhDGtuGmusGU0a2}
\end{EQA}
Moreover, with \( \rho \eqdef \dmax^{2} \dltwu_{3} \, \| \DFN \IFN^{-1} \Av \| \leq \rhos \) and
\( \rell \eqdef \mfrac{2}{1 + \sqrt{1 - 2 \rho}} \),
\begin{EQA}[rcl]
    - \mfrac{\dltwu_{3} \, \rell^{2}(3 - \rell)}{6} \| \DFN \, \IFN^{-1} \Av \|^{3}
    & \leq &  
    2 \fn(\upsvn) - 2 \fn(\upsvs) + \| \IFN^{-1/2} \Av \|^{2} 
    \leq 
    \mfrac{\dltwu_{3}}{3} \, \| \DFN \, \IFN^{-1} \Av \|^{3} 
    \, .
    \qquad
\label{3d3Af12DGttGa2}
\end{EQA}
In particular,
\begin{EQA}[c]
    \Bigl| 2 \fn(\upsvn) - 2 \fn(\upsvs) + \| \IFN^{-1/2} \Av \|^{2} \Bigr| 
    \leq 
    \mfrac{1}{2} \, \dltwu_{3} \, \| \DFN \, \IFN^{-1} \Av \|^{3} 
    \, .
\end{EQA}
\end{theorem}

\begin{proof}
W.l.o.g. 
replace \( \rr \) with \( \mfrac{4}{3} \, \dmax \| \IFN^{-1/2} \Av \| \).
\nameref{LL3tref} and Lemma~\ref{LdltwLa3t} ensure \eqref{dtb3u1DG2d3GPg} with 
\( \dltwbd(\upsv) = \dltwu_{3} \, \rr/2 \).
By \eqref{yxdhewndu7jwnjjuu}, it holds \( \dmax^{2} \dltwu_{3} \, \rr = 4 \rhos/3 \leq 1/2 \),
and \eqref{yxdhewndu7jwnjjuu} implies \eqref{rrm23r0ut3ua}.
Now the first statement follows from Theorem~\ref{Pconcgeneric}.

As \( \nabla \fs(\upsvs) = 0 \), \nameref{LL3tref} implies for any \( \upsv \in \CA(\rrn) \)
\begin{EQA}
	\bigl| 
		\fs(\upsv) - \fs(\upsvs) - \mfrac{1}{2} \| \IFN^{1/2} (\upsv - \upsvs) \|^{2} 
	\bigr|
	& \leq &
	\mfrac{\dltwu_{3}}{6} \| \DFN (\upsv - \upsvs) \|^{3}
	\, .
	\qquad \quad
\label{d3GrGELGtsG122}
\end{EQA}
Further, 
\begin{EQA}[rcl]
	&& \hspace{-12pt}
	\fn(\upsv) - \fn(\upsvs) + \mfrac{1}{2} \| \IFN^{-1/2} \Av \|^{2}
	\\
	&=&
	\bigl\langle \upsv - \upsvs, \Av \bigr\rangle
	+ \fs(\upsv) - \fs(\upsvs) + \mfrac{1}{2} \| \IFN^{-1/2} \Av \|^{2} 
	\\
	&=&
	\mfrac{1}{2} \bigl\| \IFN^{1/2} (\upsv - \upsvs) + \IFN^{-1/2} \Av \bigr\|^{2}
	+ \fs(\upsv) - \fs(\upsvs) - \mfrac{1}{2} \| \IFN^{1/2} (\upsv - \upsvs) \|^{2} .
	\qquad 
\label{12ELGuELusG2}
\end{EQA}
As \( \upsvn \in \CA(\rrn) \) and it minimizes \( \fn(\upsv) \), we derive by \eqref{d3GrGELGtsG122} and Lemma~\ref{Llin23}
with \( \AFN = \DFN^{-1} \IFN \, \DFN^{-1} \), \( \afv = - \DFN \, \IFN^{-1} \Av \), \( \rho = \dltwu_{3} \| \afv \| \),
and \( \rell \eqdef \mfrac{2}{1 + \sqrt{1 - 2 \rho}} \)
\begin{EQA}
	&& \hspace{-12pt}
	2 \fn(\upsvn) - 2 \fn(\upsvs) + \| \IFN^{-1/2} \Av \|^{2}
	=
	\inf_{\upsv \in \CA(\rrn)} 
	\Bigl\{ 
		2 \fn(\upsv) - 2 \fn(\upsvs) + \| \IFN^{-1/2} \Av \|^{2} 
	\Bigr\}
	\\
	& \geq &
	\inf_{\upsv \in \CA(\rrn)} 
	\Bigl\{ 
		\bigl\| \IFN^{1/2} (\upsv - \upsvs) + \IFN^{-1/2} \Av \bigr\|^{2} 
		- \mfrac{\dltwu_{3}}{3} \| \DFN (\upsv - \upsvs) \|^{3}
	\Bigr\} 
	\\
	& \geq &
	- \mfrac{\dltwu_{3} \, \rell^{2}(3 - \rell)}{6} \, \| \DFN \, \IFN^{-1} \Av \|^{3} \, .
\label{d3G1212222B} 
\end{EQA}
Similarly,
\begin{EQA}
	&&  
	2 \fn(\upsvn) - 2 \fn(\upsvs) + \| \IFN^{-1/2} \Av \|^{2}
	\\
	&& \quad
	\leq 
	\inf_{\upsv \in \CA(\rrn)} 
	\Bigl\{ 
		\bigl\| \IFN^{1/2} (\upsv - \upsvs) + \IFN^{-1/2} \Av \bigr\|^{2} 
		+ \mfrac{\dltwu_{3}}{3} \| \DFN (\upsv - \upsvs) \|^{3}
	\Bigr\}
	\leq 
	\mfrac{\dltwu_{3} }{3} \, \| \DFN \, \IFN^{-1} \Av \|^{3} \, . 
	\qquad \quad
\label{fd3G122B2m2} 
\end{EQA}
This completes the proof of \eqref{3d3Af12DGttGa2}.
The last statement follows 
by \( \mfrac{\rell^{2}(3 - \rell)}{6} \leq 40/81 < \mfrac{1}{2} \) for \( \rho \in [0,3/8] \)
and \( \rell \in [1,4/3] \).
\end{proof}


\Subsection{Refined bounds under locally uniform smoothness}
The bounds of Theorem~\ref{Pconcgeneric2} can be made more accurate if \( \fs \) satisfies
\nameref{LLsT3ref} and one can apply Taylor's expansion around any point close to \( \upsvs \).
In particular, the location and remainder bounds in the improved results only involve the value \( \| \DFN \IFN^{-1} \Av \| \) 
and avoid \( \| \IFN^{-1/2} \Av \| \) which can be very large in some settings; 
see Section~\ref{Slinquadr}.


\begin{theorem}
\label{PFiWigeneric2}
Let 
\begin{mylist}
\item
\( \fs(\upsv) \) be a strongly convex function with \( \fs(\upsvs) = \inf_{\upsv} \fs(\upsv) \)  
and \( \IFN = \nabla^{2} \fs(\upsvs) \);
\item
\( \fn(\upsv) \) and \( \fs(\upsv) \) be related by \eqref{4hbh8njoelvt6jwgf09} with some vector \( \Av \);
\item
the function \( \fs(\upsv) \) satisfy \nameref{LLsT3ref} at \( \upsvs \) with 
\( \DFN^{2} \), \( \rrn = \rr(\upsvs) \), and \( \dltwu_{3} \); 

\item
for some \( \dmax > 0 \) 
\begin{EQA}[c]
	\DFN^{2} \leq \dmax^{2} \hspm \IFN ,
	\quad
	\rrn \geq \mfrac{4}{3} \| \DFN \, \IFN^{-1} \Av \| \, ,
	\quad
	\dmax^{2} \hspm \dltwu_{3} \, \| \DFN \, \IFN^{-1} \Av \| \leq \mfrac{3}{8} 
	\, .
\label{8difiyfc54wrboer7bjfr}
\end{EQA}
\end{mylist}
Define
\begin{EQA}[rcl]
	\rbias
	\eqdef
	\dmax^{2} \hspm \dltwu_{3} \, \| \DFN \, \IFN^{-1} \Av \|
	\, ,
	\qquad
	\rell
	& \eqdef &
	\frac{2}{1+\sqrt{1-2\rbias}}
	\, .
\label{uyhcvdicvu998655rudvfd}
\end{EQA}
Then \( \rell \in [1,4/3] \) and the following bounds are valid:
\begin{EQA}[rcccc]
\label{dvue6554d5rdtehes}
    \| \DFN (\upsvn - \upsvs) \|
    & \leq &
    \rell \, \| \DFN \, \IFN^{-1} \Av \|
    & \leq &
    \mfrac{4}{3} \, \| \DFN \, \IFN^{-1} \Av \| 
    \, ,    
    \qquad
    \\
    \| \DFN^{-1} \IFN (\upsvn - \upsvs + \IFN^{-1} \Av) \|
    & \leq &
    \mfrac{\rell^{2}}{2} \, \dltwu_{3} \| \DFN \, \IFN^{-1} \Av \|^{2}
	& \leq &
	\mfrac{8}{9} \, \dltwu_{3} \| \DFN \, \IFN^{-1} \Av \|^{2}
	\, .
	\qquad
\label{DGttGtsGDGm13rGa2}
\end{EQA}
Also,
\begin{EQA}[rcccl]
	- \mfrac{\rell^{2} (3-\rell)}{6} \, \dltwu_{3} \| \DFN \, \IFN^{-1} \Av \|^{3}
	& \leq &
	2 \fn(\upsvn) - 2 \fn(\upsvs) + \| \IF^{-1/2} \Av \|^{2}
	& \leq &
	\mfrac{1}{3} \dltwu_{3} \| \DFN \, \IFN^{-1} \Av \|^{3}
	\, ,
	\qquad
\label{jchucyuyejr634jmdf}
\end{EQA}
yielding in particular
\begin{EQA}[rcl]
	\bigl| 2 \fn(\upsvn) - 2 \fn(\upsvs) + \| \IF^{-1/2} \Av \|^{2} \bigr|
	& \leq &
	\mfrac{1}{2} \, \dltwu_{3} \| \DFN \, \IFN^{-1} \Av \|^{3}
	\, .
\end{EQA}
\end{theorem}

\begin{remark}[Sharpness]
For the scalar local cubic \( \fs(x)=x^{2}/2 - \dltwu_{3} x^{3}/6 \), 
\( \DFN = \IFN = \dmax = 1 \), let \( A > 0 \) and \(\dltwu_{3} A \uparrow 3/8 \). 
The minimizer \( \xd \) of \( \fn(x) = \fs(x) + A x \) satisfies \( \xd/A \uparrow 4/3 \). 
The ratio of the location remainder to \(\dltwu_{3} A^{2} \) tends to \( 8/9 \), 
and the ratio of the absolute doubled-value remainder to \( \dltwu_{3} A^{3} \) tends to \( \rell^{2} (3 - \rell)/6 = 40/81 \)
for \( \rell = 4/3 \). 
Consequently, uniform constants \( \rell \leq \mfrac{4}{3} \), \( \mfrac{\rell^{2}}{2} \leq \mfrac{8}{9} \),
\( \mfrac{\rell^{2} (3-\rell)}{6} \leq \frac{40}{81} \) and \( \mfrac{1}{3} \) 
in \eqref{dvue6554d5rdtehes} through \eqref{jchucyuyejr634jmdf} 
cannot be lowered under the theorem's assumptions.
\end{remark}

\begin{proof}
If the function \( \fs \) is quadratic and convex with the minimum at \( \upsvs \) then the linearly perturbed function
\( \fn \) is also quadratic and convex with the minimum at \( \upsvf = \upsvs - \IFN^{-1} \Av \).
In general, the point \( \upsvf \) is not the minimizer of \( \fn \), however, it is very close to \( \upsvn \).
We use that \( \nabla \fs(\upsvs) = 0 \) 
and \( \nabla^{2} \fs(\upsvs) = \IFN \).
The definition \eqref{uyhcvdicvu998655rudvfd} of \( \rell \) and \( \rbias \in [0,3/8] \) imply
\begin{EQA}[c]
\label{eq:p34-ell}
 	1 \leq \rell \leq \mfrac{4}{3},
 	\qquad 
	\rell - 1 = \mfrac{\rbias \, \rell^{2}}{2}
	\, ,
 	\qquad 
	\rbias \, \rell
	=
	1 - \sqrt{1-2\rbias}
	\leq 
	\mfrac{1}{2}
	\, .
\label{jviuu8767jfloey7eufg}
\end{EQA}
These identities also hold at \( \rbias=0 \), where \( \rell=1 \).

For \( \BBH \eqdef \DFN^{-1} \IFN \DFN^{-1} \), 
it holds \(\| \BBH^{-1} \| \leq \dmax^{2} \).
Define \( \afv \eqdef - \BBH^{-1} \DFN^{-1} \Av = - \DFN \, \IFN^{-1} \Av \)
and introduce \( \uv \eqdef \DFN (\upsv - \upsvs) \) and
\begin{EQA}[c]
	R(\uv)
	\eqdef
	\DFN^{-1} \{ \nabla \fs(\upsv) - \IF (\upsv - \upsvs) \}
	=
	\DFN^{-1} \nabla \fs(\upsvs + \DFN^{-1} \uv) - \BBH \uv
	\, .
\end{EQA}
By \eqref{y6sdjsdy7erwmcuecuid} and \eqref{jhfuy7f7dfyedye663eh} of Lemma~\ref{Ldltw4s} ahead,
\begin{EQA}[c]
	\| R(\uv) \|
	\leq
	\mfrac{1}{2} \, \dltwu_{3} \, \| \DFN (\upsv - \upsvs) \|^{2}
	=
	\mfrac{1}{2} \, \dltwu_{3} \, \| \uv \|^{2}
	\, ,
	\qquad 
	\| \nabla R(\uv) \|
	\leq 
	\dltwu_{3} \, \| \uv \|
	\, .
\label{huy77ewjdf766hkgf}
\end{EQA}
The stationary equation \( \nabla \fn(\upsvn) = 0 \) leads to \( \nabla \fs(\upsvn) = - \Av \) or
\begin{EQA}[rcl]
	\uvn
	\eqdef
	\DFN (\upsvn - \upsvs) 
	&=& 
	- \DFN \, \IFN^{-1} \Av - \DFN \, \IFN^{-1} \DFN \, R(\uvn) 
	=
	\afv - \BBH^{-1} R(\uvn)
	\, .
\end{EQA} 
Define \( \Phi(\uv) \eqdef \afv - \BBH^{-1} R(\uv) \).
Then the latter identity is the fixed-point equation \( \uvn = \Phi(\uvn) \).
On the ball of radius \( \rrd = \rell \, \| \afv \| \), by \( \| \BBH^{-1} \| \leq \dmax^{2} \), 
\eqref{uyhcvdicvu998655rudvfd}, and \eqref{huy77ewjdf766hkgf},
\begin{EQA}[c]
	\| \Phi(\uv) \|
	\leq 
	\| \DFN \, \IFN^{-1} \Av \| + \| \DFN \, \IFN^{-1} \DFN \| \; \| R(\uv) \|
	\leq 
	\| \afv \| + \mfrac{\dmax^{2} \dltwu_{3}}{2} \rrd^{2}
	=
	\| \afv \| \bigl( 1 + \mfrac{\rho \rell^{2}}{2} \bigr)
	=
	\rrd
	\, ,
	\\
	\| \nabla \Phi(\uv) \|
	\leq 
	\dmax^{2} \dltwu_{3} \, \rrd
	=
	\rbias \, \rell
	\leq 
	1/2
	\, .
\end{EQA}
Therefore, the map \( \Phi \) has a unique fixed point in the ball. 
Strong convexity ensures that the stationary point is the global minimizer. 
This proves \eqref{dvue6554d5rdtehes} by \( \rell \leq 4/3 \).

Further, stationarity and \eqref{y6sdjsdy7erwmcuecuid} give
\begin{EQA}[rcl]
	\| \DFN^{-1} \IFN (\upsvn - \upsvs + \IFN^{-1} \Av) \|
	&=&
	\| R(\uvn) \|
	\leq
	\mfrac{1}{2} \, \dltwu_{3} \rrd^{2}
	=
	\mfrac{1}{2} \, \dltwu_{3} \rell^{2} \| \afv \|^{2}
	\, ,
\end{EQA}
and \eqref{DGttGtsGDGm13rGa2} follows.
The bounds in \eqref{jchucyuyejr634jmdf} are proved as \eqref{3d3Af12DGttGa2} in Theorem~\ref{Pconcgeneric2}.
\end{proof}

\begin{lemma}
\label{Ldltw4s}
Assume \nameref{LLsT3ref} at \( \upsv \)
and let \( \UVz_{\rr} = \{ \uv \colon \| \DFN \uv \| \leq \rr \} \).
Then 
\begin{EQA}[ccl]
	\bigl\| \DFN^{-1} \bigl\{ \nabla \fs(\upsv + \uv) - \nabla \fs(\upsv) 
	- \nabla^{2} \fs(\upsv) \uv  \bigr\} 
	\bigr\|
	& \leq &
	\mfrac{1}{2} \, \dltwu_{3} \, \| \DFN \uv \|^{2} \, ,
	\quad
	\uv \in \UVz_{\rr} \, .
	\qquad
\label{y6sdjsdy7erwmcuecuid}
\end{EQA}
Also for all \( \uv, \uv_{1} \in \UVz_{\rr} \) 
\begin{EQ}[rcl]
	\bigl\| \DFN^{-1} \bigl\{ \nabla^{2} \fs(\upsv + \uv_{1}) - \nabla^{2} \fs(\upsv + \uv) \bigr\} \DFN^{-1}	\bigr\|
	& \leq &
	\dltwu_{3} \, \| \DFN (\uv_{1} - \uv) \|
	\, ,
\label{jhfuy7f7dfyedye663eh}
\end{EQ}
and
\begin{EQA}
	&& \nquad
	\bigl\| \DFN^{-1} \bigl\{ \nabla \fs(\upsv + \uv_{1}) - \nabla \fs(\upsv + \uv) - \nabla^{2} \fs(\upsv) (\uv_{1} - \uv) \bigr\}
	\bigr\|
	\\
	& \leq &
	\dltwu_{3} \| \DFN (\uv_{1} - \uv) \| \bigl( \| \DFN \uv \| + \mfrac{1}{2} \| \DFN (\uv_{1} - \uv) \| \bigr) 
	\, .
	\qquad
\label{6dcfujcu8ed8edsudyf5tre35}
\end{EQA}
Moreover, under \nameref{LLsT4ref}, for any \( \uv \in \UVz_{\rr} \),
\begin{EQA}
	\bigl\| \DFN^{-1} \bigl\{ \nabla \fs(\upsv + \uv) - \nabla \fs(\upsv) 
	- \nabla^{2} \fs(\upsv) \uv  
	- \mfrac{1}{2} \langle \nabla^{3} \fs(\upsv), \uv^{\otimes 2} \rangle \bigr\} \bigr\|
	& \leq &
	\mfrac{1}{6} \dltwu_{4} \, \| \DFN \uv \|^{3} \, .
	\qquad
\label{y6sdjsdy7erwmcuecuid4}
\end{EQA}
\end{lemma}

\begin{proof}
Fix \( \uv \in \UVz_{\rr} \) and any vector \( \wv \in \R^{\dimp} \).
For \( t \in [0,1] \), denote
\begin{EQA}
	h(t)
	& \eqdef &
	\bigl\langle \nabla \fs(\upsv + t\uv) - \nabla \fs(\upsv) 
	- t \langle \nabla^{2} \fs(\upsv), \uv \rangle, \, \wv \bigl\rangle
	\, .
\label{7dcjef6chjfer7v54etghf}
\end{EQA}
Then \( h(0) = 0 \), \( h'(0) = 0 \), and 
\nameref{LLsT3ref} and \eqref{jbuyfg773jgion94euyyfg} imply 
\begin{EQA}[rcl]
	|h''(t)|
	&=&
	\bigl| \langle \nabla^{3} \fs(\upsv + t \uv), \uv^{\otimes 2} \otimes \wv \rangle \bigr|
	\leq 
	\dltwu_{3} \| \DFN \uv \|^{2} \, \| \DFN \wv \|
	\, .
\end{EQA}
With \( \av \eqdef \nabla \fs(\upsv + \uv) - \nabla \fs(\upsv) 
	- \nabla^{2} \fs(\upsv) \uv \), this yields
\begin{EQA}
	|\langle \av,\wv \rangle| 
	&=& 
	|h(1)| \leq \mfrac{1}{2} \, \dltwu_{3} \| \DFN \uv \|^{2} \, \| \DFN \wv \|
	\, ,
	\\
	\| \DFN^{-1} \av \|
	&=&
	\sup_{\| \wv \| = 1} \bigl| \langle \DFN^{-1} \av,\wv \rangle \bigr|
	=
	\sup_{\| \wv \| = 1} \bigl| \langle \av,\DFN^{-1} \wv \rangle \bigr|
	\leq 
	\mfrac{1}{2} \, \dltwu_{3} \, \| \DFN \uv \|^{2} 
	\, ,
\label{udfjd6vhfwe36vneo}
\end{EQA}
and the first statement follows. 
To prove \eqref{y6sdjsdy7erwmcuecuid4}, define
\begin{EQA}
	\av
	& \eqdef &
	\nabla \fs(\upsv + \uv) - \nabla \fs(\upsv) 
	- \nabla^{2} \fs(\upsv) \uv  
	- \mfrac{1}{2} \langle \nabla^{3} \fs(\upsv), \uv^{\otimes 2} \rangle
	\, ,
	\\
	h(t)
	& \eqdef &
	\bigl\langle \nabla \fs(\upsv + t\uv) - \nabla \fs(\upsv) 
	- t \nabla^{2} \fs(\upsv) \uv  
	- \mfrac{t^{2}}{2} \langle \nabla^{3} \fs(\upsv), \uv^{\otimes 2} \rangle, \wv \bigl\rangle 
	\, ,
\label{t6dtwsghwesyfyghe322w2w}
\end{EQA} 
and apply 
\eqref{jbuyfg773jgion94euyyfg4}.

Further, with \( \Delta = \uv_{1} - \uv \), 
by \nameref{LLsT3ref}, for any \( \wv \in \R^{\dimp} \) and some \( t \in [0,1] \),
\begin{EQA}
	&& \nquad
	\bigl| \langle \DFN^{-1} \bigl\{ \nabla^{2} \fs(\upsv + \uv_{1}) - \nabla^{2} \fs(\upsv + \uv) \bigr\} \, \DFN^{-1}, \wv^{\otimes 2} \rangle \bigr|
	\\
	& = &
	\bigl| \bigl\langle \nabla^{3} \fs(\upsv + \uv + t \Delta), \Delta \otimes (\DFN^{-1} \wv)^{\otimes 2} \bigr\rangle \bigr|
	\leq 
	\dltwu_{3} \| \DFN \Delta \| \, \| \wv \|^{2}
	\, .
\label{d76jnwef8j2egtftftyjr5}
\end{EQA}
This proves \eqref{jhfuy7f7dfyedye663eh}.

Further, 
for any \( \uv,\uv_{1} \in \UVz_{\rr} \), it holds with \( \Delta \eqdef \uv_{1} - \uv \), 
\begin{EQA}
 	&& \nquad
	\DFN^{-1} \bigl\{ \nabla \fs(\upsv + \uv_{1}) - \nabla \fs(\upsv + \uv) - \nabla^{2} \fs(\upsv) \Delta \bigr\} 
	\\
	& = &
	\int_{0}^{1} \DFN^{-1} \bigl\{ \nabla^{2} \fs(\upsv + \uv + t \Delta) - \nabla^{2} \fs(\upsv) \bigr\} \DFN^{-1} \, \DFN \Delta \; dt
	\, ,
\end{EQA}
and by \eqref{jhfuy7f7dfyedye663eh}, 
\begin{EQA}[rcl]
 	&& \nquad
	\bigl\| \DFN^{-1} \bigl\{ \nabla \fs(\upsv + \uv_{1}) - \nabla \fs(\upsv + \uv) - \nabla^{2} \fs(\upsv) \Delta \bigr\} \bigr\| 
	\\
	& \leq &
 	\dltwu_{3} \| \DFN \Delta \| \int_{0}^{1} \| \DFN \uv + t \DFN \Delta \| \, dt
 	\leq 
 	\dltwu_{3} \| \DFN \Delta \| (\| \DFN \uv \| + \mfrac{1}{2} \| \DFN \Delta \|) 
	\, ,
\end{EQA}
and \eqref{6dcfujcu8ed8edsudyf5tre35} follows.
\end{proof}


\Subsection{{Fourth-order expansions}}
Under fourth-order condition \nameref{LLsT4ref},  expansion \eqref{DGttGtsGDGm13rGa2} can be further refined.  

\begin{theorem}
\label{Pconcgeneric4}
Let 
\begin{mylist}
\item
\( \fs(\upsv) \) be a strongly convex function with \( \fs(\upsvs) = \inf_{\upsv} \fs(\upsv) \)  
and \( \IFN = \nabla^{2} \fs(\upsvs) \);
\item
\( \fn(\upsv) \) and \( \fs(\upsv) \) be related by \eqref{4hbh8njoelvt6jwgf09} with some vector \( \Av \);
\item
\( \fs(\upsv) \) satisfy \nameref{LLsT3ref} and \nameref{LLsT4ref} at \( \upsvs \) with 
\( \DFN^{2} \), \( \dltwu_{3} \), \( \dltwu_{4} \), and \( \rrn = \rr(\upsvs) \); 

\item
for some \( \dmax > 0 \) 
\begin{EQA}[c]
	\DFN^{2} \leq \dmax^{2} \hspm \IFN 
	\, ,
	\quad
	\rrn \geq \mfrac{4}{3} \| \DFN \, \IFN^{-1} \Av \| 
	\, ,
	\quad
	\dmax^{2} \hspm \dltwu_{3} \| \DFN \, \IFN^{-1} \Av \| \leq \mfrac{3}{8} 
	\, .
	\qquad
\label{8difiyfc54wrbosT4}
\end{EQA}
\end{mylist}
Then \( \| \DFN (\upsvn - \upsvs) \| \leq (4/3) \| \DFN \, \IFN^{-1} \Av \| \). 
Further, define
\begin{EQA}
	\avn 
	&=&
	- \IFN^{-1} \{ \Av + \nabla \Tens(\IFN^{-1} \Av) \} \, ,
\label{8vfjvr43223efryfuweef}
\end{EQA}
where \( \Tens(\uv) = \frac{1}{6} \langle \nabla^{3} \fs(\upsvs), \uv^{\otimes 3} \rangle \) for \( \uv \in \R^{\dimp} \).
Then with \( \rell \) from \eqref{uyhcvdicvu998655rudvfd}
\begin{EQA}
  	\| \DFN^{-1} \IFN (\upsvn - \upsvs - \avn) \|
  	& \leq &
  	\Bigl\{ \mfrac{1}{6} \dltwu_{4} + \mfrac{\rell^{2}(\rell+1)}{4} \, \dmax^{2} \dltwu_{3}^{2} \Bigr\} \| \DFN \, \IFN^{-1} \Av \|^{3}
  	\\
	& \leq &
  	\Bigl( \mfrac{1}{6} \dltwu_{4} + \mfrac{28}{27} \dmax^{2} \dltwu_{3}^{2} \Bigr) \| \DFN \, \IFN^{-1} \Av \|^{3}
    \, .
 	\qquad
\label{DGttGtsGDGm13rGa4}
\end{EQA}
For the value remainder
\begin{EQA}[rcl]
	\Reminder
	& \eqdef &
	\fn(\upsvn) - \fn(\upsvs) + \frac{1}{2} \| \IFN^{-1/2} \Av \|^{2} + \Tens(\IFN^{-1} \Av)
	\, ,
\label{jhu89erklvytr68urtgr}
\end{EQA}
it holds 
\begin{EQA}[c]
\label{3d3Af12DGttGa4}
    - \left\{ \mfrac{1}{24} \dltwu_{4} + C(\rell) \dmax^{2} \dltwu_{3}^{2} \right\} \| \DFN \, \IFN^{-1} \Av \|^{4}
    \leq 
    \Reminder 
    \leq 
    \mfrac{1}{24} \, \dltwu_{4} \, \| \DFN \, \IFN^{-1} \Av \|^{4}
    \, ,
    \\
    C(\rell)
    \eqdef 
    \mfrac{1}{24} \, \rell^{2}(2+2\rell-\rell^{2})
    \, .
\label{3d3Af12DGttGa4T}
\end{EQA}
In particular,
\begin{EQA}[c]
 	|\Reminder|
	\leq
	\Bigl( \mfrac{1}{24} \dltwu_{4} + \mfrac{52}{243} \dmax^{2} \dltwu_{3}^{2} \Bigr) \| \DFN \, \IFN^{-1} \Av \|^{4}.
\label{3d3Af12DGttGa4eee}
\end{EQA}
Also 
\begin{EQA}
	\bigl| \Tens(\IFN^{-1} \Av) \bigr|
	& \leq &
	\mfrac{1}{6} \, \dltwu_{3} \| \DFN \, \IFN^{-1} \Av \|^{3} 
	\, .
\label{u7jcc7e45hfiobvioeye6yhy}
\end{EQA}
\end{theorem}

\begin{proof}
If \( \biasl \eqdef \| \DFN \, \IFN^{-1} \Av \| = 0 \), then \( \Av = 0 \), and all assertions follow immediately.
Assume \( \biasl > 0 \). By \eqref{uyhcvdicvu998655rudvfd} and \eqref{jviuu8767jfloey7eufg}
\begin{EQA}[c]
	1 \leq \rell \leq \mfrac{4}{3}
	\, ,
	\qquad
	\rell-1 = \mfrac{1}{2} \, \dmax^{2} \dltwu_{3} \, \rell^{2} \biasl
	\, ,
	\qquad
	\dmax^{2} \dltwu_{3} \, \rell \, \biasl \leq \mfrac{1}{2}
	\, .
\label{uhcdi89s977ufosewdf}
\end{EQA}
Theorem~\ref{PFiWigeneric2} gives with \( \biasl = \| \DFN \, \IFN^{-1} \Av \| \)
\begin{EQA}[c]
	\| \DFN (\upsvn-\upsvs) \|
 	\leq 
	\rell \, \biasl 
	\, ,
	\qquad
	\| \DFN^{-1} \IFN (\upsvn - \upsvs + \IFN^{-1} \Av) \|
    \leq 
    \mfrac{1}{2} \, \dltwu_{3} \, \rell^{2} \biasl^{2}
	\, .
\label{ihgfe9iythjkuuf}
\end{EQA}
Denote 
\begin{EQA}[rcl]
	\upsvf
	& \eqdef &
	\upsvs - \IF^{-1} \Av
	\, .
\end{EQA}
For \( \BFN = \DFN^{-1} \IFN \DFN^{-1} \), 
the condition \( \| \BFN^{-1} \| = \| \DFN \IFN^{-1} \DFN \| \leq \dmax^{2} \) and \eqref{uhcdi89s977ufosewdf} yield
\begin{EQA}[rcl]
	\| \DFN (\upsvn - \upsvf) \|
	&=&
	\| \DFN \IFN^{-1} \DFN \; \DFN^{-1} \IFN (\upsvn - \upsvf) \|
	\leq 
	\| \BFN^{-1} \| \; \| \DFN^{-1} \IFN (\upsvn - \upsvf) \|
	\\
	& \leq &
	\dmax^{2} \|\DFN^{-1} \IFN (\upsvn - \upsvs + \IFN^{-1}\Av)\|
	\leq
	\mfrac{1}{2} \, \dmax^{2} \dltwu_{3} \, \rell^{2} \biasl^{2}
 	= 
	(\rell-1) \biasl
	\, .
	\qquad
\label{iiiiikdeyfhjerjkfer}
\end{EQA}
Thus, \( \upsvf \), \( \upsvn \), and the segment joining
them belong to the neighborhood on which the smoothness
assumptions hold.

Later we use that \( \Tens(- \uv) = - \Tens(\uv) \) while \( \nabla \Tens(-\uv) = \nabla \Tens(\uv) \).
By \eqref{y6sdjsdy7erwmcuecuid4} of Lemma~\ref{Ldltw4s} and \( \nabla \fs(\upsvs) = 0 \)
\begin{EQA}[rcl]
	\bigl\| \DFN^{-1} \bigl\{ \nabla \fs(\upsvf) + \Av - \nabla \Tens(\IFN^{-1} \Av) \bigr\} \bigr\|
	& \leq &
	\mfrac{1}{6} \dltwu_{4} \| \DFN \IFN^{-1} \Av \|^{3}
	=
	\mfrac{1}{6} \dltwu_{4} \, \biasl^{3}
	\, .
\label{jvgi8uytrhgyh3ekbueh}
\end{EQA}
By \eqref{6dcfujcu8ed8edsudyf5tre35} and \eqref{uhcdi89s977ufosewdf}
in view of \( \| \DFN (\upsvf - \upsvs) \| = \| \DFN \IFN^{-1} \Av \| = \biasl \), it holds by \eqref{iiiiikdeyfhjerjkfer}
\begin{EQA}[rcl]
	&& \nquad
	\bigl\| \DFN^{-1} \bigl\{ \nabla \fs(\upsvf) - \nabla \fs(\upsvd) - \IFN (\upsvf - \upsvd) \bigr\}
	\bigr\|
	\\
	& \leq &
	\dltwu_{3} \bigl( \| \DFN (\upsvf - \upsvs) \| + \mfrac{1}{2} \| \DFN (\upsvd - \upsvf) \| \bigr) \; \| \DFN (\upsvd - \upsvf) \| 
	\\
	& \leq &
	\dltwu_{3} \Bigl\{ \biasl + \mfrac{1}{2} (\rell-1) \biasl \Bigr\} \mfrac{1}{2} \, \dmax^{2} \dltwu_{3} \, \rell^{2} \biasl^{2}
	=
	\mfrac{1}{4} \, \rell^{2} (\rell + 1) \dmax^{2} \dltwu_{3}^{2} \, \biasl^{3}
	\, .
	\qquad
\end{EQA}
Further, \( \nabla \fn(\upsvd) = 0 \) implies \( \nabla \fs(\upsvd) + \Av = 0 \).
This and \eqref{jvgi8uytrhgyh3ekbueh} prove \eqref{DGttGtsGDGm13rGa4}, using \( \rell^{2} (\rell+1) / 4 \leq 28/27 \).

We next prove the value bounds.
Taylor expansion at \( \upsvs \), together with \nameref{LLsT4ref}, gives
\begin{EQA}[c]
	\bigl| \fn(\upsvf)-\fn(\upsvs) + \mfrac{1}{2} \|\IFN^{-1/2}\Av\|^{2} + \Tens(\IFN^{-1}\Av) \bigr|
 	\leq 
	\mfrac{1}{24} \, \dltwu_{4} \, \biasl^{4}
	\, .
\end{EQA}
Here we used
\( \nabla \fn(\upsvs) = \Av \),
\( \nabla^{2} \fn(\upsvs) = \IFN \), and
\( \Tens(-\IFN^{-1} \Av) = -\Tens(\IFN^{-1} \Av) \).
Since \( \fn(\upsvn) \leq \fn(\upsvf) \), this implies
\begin{EQA}[c]
	\Reminder \leq \mfrac{1}{24} \, \dltwu_{4} \, \biasl^{4}
	\, .
\end{EQA}
For the lower bound, we estimate
\( \fn(\upsvf)-\fn(\upsvn) \).
By \eqref{y6sdjsdy7erwmcuecuid} with \( \uv = - \IFN^{-1} \Av \),
\begin{EQA}[c]
	\| \DFN^{-1} \nabla \fn(\upsvf) \|
 	= 
	\bigl\| \DFN^{-1} \bigl\{ \nabla \fs(\upsvf) + \Av \bigr\} \bigr\|
 	\leq 
	\mfrac{1}{2} \, \dltwu_{3} \, \| \DFN \IFN^{-1} \Av \|^{2}
	=
	\mfrac{1}{2} \, \dltwu_{3} \, \biasl^{2}
	\, .
\end{EQA}
Furthermore, for \( 0 \leq t \leq 1 \),
\begin{EQA}[c]
	\bigl\| \DFN\bigl\{ \upsvf + t (\upsvn-\upsvf) - \upsvs \bigr\} \bigr\|
	\leq
	\biasl + t \| \DFN (\upsvn-\upsvf) \|
 	\leq 
	\rell \, \biasl
	\leq 
	\rrn
	\, .
\end{EQA}
Consequently, \eqref{jhfuy7f7dfyedye663eh} implies
\begin{EQA}[rl]
	\nabla^{2} \fn\bigl(\upsvf + t (\upsvn - \upsvf) \bigr)
	\succeq
	\IFN - \dltwu_{3} \bigl( \biasl + t \| \DFN (\upsvn-\upsvf) \| \bigr) \DFN^{2}
	\, .
\end{EQA}
Taylor's formula with integral remainder therefore gives by \( \IFN\succeq\dmax^{-2}\DFN^{2} \)
\begin{EQA}[rcl]
	&& \nquad
	\fn(\upsvn)-\fn(\upsvf)
	\\
	& = &
	\bigl\langle \nabla \fn(\upsvf) , \upsvn-\upsvf \bigr\rangle
	+ \int_{0}^{1}(1-t) 
	\bigl\langle \nabla^{2}\fn\bigl( \upsvf+t(\upsvn-\upsvf) \bigr), (\upsvn-\upsvf)^{\otimes2} \bigr\rangle \, dt
	\\
	& \geq &
	- \frac{\dltwu_{3}\biasl^{2}}{2} \| \DFN (\upsvn-\upsvf) \|
	+ \frac{\dmax^{-2}-\dltwu_{3}\biasl}{2} \| \DFN (\upsvn-\upsvf) \|^{2}
	- \frac{\dltwu_{3}}{6} \, \| \DFN (\upsvn-\upsvf) \|^{3}
	\, .
	\qquad
\label{jucviui8vc7h3f6f6r}
\end{EQA}
Introduce the polynomial
\begin{EQA}[c]
	P(t)
	\eqdef
	\frac{\dltwu_{3}\biasl^{2}}{2} \, t
	- \frac{\dmax^{-2} - \dltwu_{3} \biasl}{2} \, t^{2}
	+ \frac{\dltwu_{3}}{6} \, t^{3}
	\, .
\end{EQA}
With this notation, \eqref{jucviui8vc7h3f6f6r} becomes
\begin{EQA}[c]
	\fn(\upsvf) - \fn(\upsvn)
	\leq 
	P\bigl( \| \DFN (\upsvn-\upsvf) \| \bigr)
	\, .
\label{jucviui8vc7h3f6f6rr}
\end{EQA}
We already know that \( \| \DFN (\upsvn-\upsvf) \| \leq (\rell-1) \biasl \).
Therefore, it suffices to bound \( P(t) \) from above on the interval \( \II \eqdef [0 , (\rell-1)\biasl] \). 
It holds on \( \II \) by
\( \rell-1 = \dmax^{2} \dltwu_{3} \, \rell^{2} \biasl/2 \); see \eqref{uhcdi89s977ufosewdf}:
\begin{EQA}[c]
	P''(t)
	=
	\dltwu_{3} (\biasl+t) - \dmax^{-2}
 	\leq 
	\dltwu_{3} \, \rell \, \biasl - \dmax^{-2}
	\leq
	- \frac{1}{2\dmax^{2}}
	\, .
\end{EQA}
Moreover,
\begin{EQA}[c]
	P'(t)
	=
	\mfrac{1}{2} \, \dltwu_{3}(\biasl + t)^{2} - \frac{t}{\dmax^{2}}
	\, ,
\end{EQA}
and \( P'(t) = 0 \) for \( t = (\rell-1) \biasl \).
This implies \( P'(t) \geq 0 \) for \( t \in \II \) and thus, 
\( P(t) \leq P\bigl( (\rell-1) \biasl \bigr) \) for \( t \in \II \);
see Lemma~\ref{L3poly}.
By \eqref{jucviui8vc7h3f6f6rr} 
\begin{EQA}[rcl]
	0 
	& \leq & 
	\fn(\upsvf)-\fn(\upsvn)
	\leq 
	P\bigl( \| \DFN (\upsvn-\upsvf) \| \bigr)
	\leq 
	P( (\rell - 1) \biasl )
	\\
	&=&
	\frac{\rell^{2}(2+2\rell-\rell^{2})}{24} \, \dmax^{2} \dltwu_{3}^{2} \, \biasl^{4}
	=
	C(\rell)\dmax^{2}\dltwu_{3}^{2} \, \biasl^{4}
	\, ,
\end{EQA}
where the equality again follows from \( \rell-1 = \dmax^{2} \dltwu_{3} \, \rell^{2} \biasl/2 \).

Combining this estimate with the Taylor bound at
\( \upsvf \) and the identity
\begin{EQA}[rcl]
	\Reminder
	& = &
	\fn(\upsvf)-\fn(\upsvs) + \frac{1}{2} \|\IFN^{-1/2}\Av\|^{2} + \Tens(\IFN^{-1}\Av) 
	- \bigl\{ \fn(\upsvf) - \fn(\upsvn) \bigr\}
\end{EQA}
proves \eqref{3d3Af12DGttGa4}:
\begin{EQA}[c]
	- \Bigl\{ \mfrac{1}{24} \, \dltwu_{4} + C(\rell) \dmax^{2}\dltwu_{3}^{2} \Bigr\} \biasl^{4}
 	\leq 
	\Reminder
 	\leq 
	\mfrac{1}{24} \, \dltwu_{4} \, \biasl^{4}
	\, .
\end{EQA}
Finally,
\begin{EQA}[c]
	C'(\rell)
 	= 
	\frac{\rell(2-\rell)(2\rell+1)}{12}
	\geq
	0
	\, ,
	\qquad 
	1 \leq \rell \leq \frac{4}{3}
	\, ,
\end{EQA}
so that \( C(\rell)\le C(4/3) = \mfrac{52}{243} \).
This proves \eqref{3d3Af12DGttGa4eee}.

The remaining assertion follows directly from \nameref{LLsT3ref} at \( \upsvs \):
\begin{EQA}[c]
	|\Tens(\IFN^{-1}\Av)|
 	= 
	\mfrac{1}{6}
	\bigl| \bigl\langle \nabla^{3}\fs(\upsvs),(\IFN^{-1}\Av)^{\otimes3} \bigr\rangle \bigr|
	\le
	\mfrac{1}{6} \, \dltwu_{3}\|\DFN \IFN^{-1}\Av\|^{3}
	\, .
\end{EQA}
This completes the proof.
\end{proof}

\begin{remark}[Sharpness]
\label{Rsharpness4}
Each coefficient in the uniform location and absolute value bounds is
optimal when the other coefficient is held fixed. To verify the
coefficients of \( \dmax^{2} \dltwu_{3}^{2} \), take
\[
 	\fs(\ups) = \frac{\ups^{2}}{2} - \frac{\dltwu \ups^{3}}{6}
	\, ,
	\qquad
 	\fn(\ups) = f(\ups) - A \ups,
	\qquad 
	\DFN = \IFN = \kappa = 1,
\]
where \( \dltwu > 0 \), \( 0 < \dltwu \, A \leq 3/8\), and
\( \Upsilon = (-3/(4\dltwu),3/(4\dltwu)) \).
Then \( \fs''(\ups) > 1/4 \), and the assumptions hold with
\( \rrn = 4A/3 \), \( \dltwu_{3} = \dltwu \), and \( \dltwu_{4} = 0 \).
Writing \( \rho = \dltwu \, A \), one has \( \biasl = A \),
\[
 	\upsn = \rell \, A,
	\qquad 
	a = A +\frac{\dltwu A^{2}}{2}
	\, ,
	\qquad
 	\frac{|\upsn - a|}{\dltwu^{2} A^{3}} = \frac{\rell^{2}(\rell+1)}{4}
	\, ,
 	\qquad 
 	\Reminder = - C(\rell) \dltwu^{2} A^{4}
	\, .
\]
At \( \rho = 3/8 \), these ratios give \(28/27\) and \(52/243\).
The same optimal constants apply under the strict condition \(\rho<3/8\)
by taking the limit \( \rho \uparrow 3/8 \).

For the coefficients of \( \dltwu_{4} \), take
\[
 	\fs(\ups) = \frac{\ups^{2}}{2} + \frac{\gamma \ups^{4}}{24},
	\qquad
 	\fn(\ups) = \fs(\ups) - A \ups,
	\qquad 
	\DFN = \IFN = \kappa = 1,
\]
with fixed \( \gamma > 0 \), \( \rrn = 4A/3 \), \( \dltwu_{4} = \gamma \), and
\( \dltwu_{3} = \gamma \rrn \). 
As \( A \downarrow 0\),
\[
 	a = b,
	\qquad 
	|\upsn - a| = \frac{\gamma}{6} A^{3} + O(A^{5})
	\, ,
 	\qquad 
	\Reminder = \frac{\gamma}{24} A^{4} + O(A^{6}).
\]
The terms involving \( \dltwu_{3}^{2} \) have higher order, so the coefficients
\( 1/6 \) and \( 1/24 \) cannot be decreased.
\end{remark}

\Section{A general perturbation}
This section explains how the results obtained for a linear perturbation can be extended a more general perturbation.
We study separately the cases of additive quadratic or smooth penalization, and a general perturbation.

\Subsection{Quadratic penalization}
\label{Slinquadr}

Here we discuss the case when \( \fn(\upsv) - \fs(\upsv) \) is quadratic.
More precisely, we limit ourselves to the situation with \( \fn(\upsv) = \fs(\upsv) + \| \GP \upsv \|^{2}/2 \)
for a symmetric semidefinite matrix \( \GP \).
To make the dependence on \( \GP \) explicit, denote 
\( \fG(\upsv) \eqdef \fs(\upsv) + \| \GP \upsv \|^{2}/2 \),
\begin{EQA}
	\upsvs 
	&=& 
	\arginf_{\upsv} \fs(\upsv),
	\quad
	\upsvs_{\GP} = \arginf_{\upsv} \fG(\upsv) 
	=
	\arginf_{\upsv} \bigl\{ \fs(\upsv) + \| \GP \upsv \|^{2}/2 \bigr\}.
	\qquad
\label{8cvkfc9fujf6jnmcer4cd}
\end{EQA}
We study the displacement \( \upsvs_{\GP} - \upsvs \) 
and the decrease \( \fG(\upsvs_{\GP}) - \fG(\upsvs) \) induced by this penalization.
The quadratic case provides exact formulas and motivates the general expansion.

\begin{lemma}
\label{Lbiasquadgen}
Let \( \fs(\upsv) \) be quadratic with \( \IFN \equiv \nabla^{2} \fs(\upsv) \).
Denote \( \IFN_{\GP} = \IFN + \GP^{2} \).
Then \( \upsvs_{\GP} \) from \eqref{8cvkfc9fujf6jnmcer4cd} satisfies
\begin{EQA}[rcccl]
	\upsvs_{\GP} - \upsvs
	&=&
	- \IFN_{\GP}^{-1} \GP^{2} \upsvs,
	\qquad
	\fG(\upsvs_{\GP}) - \fG(\upsvs)
	&=&
	- \frac{1}{2} \| \IFN_{\GP}^{-1/2} \GP^{2} \upsvs \|^{2} \, .
\label{kjfydf554dfwertdgwdf}
\end{EQA}
\end{lemma} 

\begin{proof}
By definition, \( \fG(\upsv) \) is quadratic with \( \nabla^{2} \fG(\upsv) \equiv \IFN_{\GP} \) and
\begin{EQA}
	\nabla \fG(\upsvs_{\GP}) - \nabla \fG(\upsvs)
	&=&
	\IFN_{\GP} \, (\upsvs_{\GP} - \upsvs) .
\label{dcudydye67e6dy3wujhdqu}
\end{EQA}
Further, \( \nabla \fs(\upsvs) = 0 \) yielding \( \nabla \fG(\upsvs) = \GP^{2} \upsvs \).
Together with \( \nabla \fG(\upsvs_{\GP}) = 0 \), this implies
\( \upsvs_{\GP} - \upsvs = - \IFN_{\GP}^{-1} \GP^{2} \upsvs \).
The Taylor expansion of \( \fG \) at \( \upsvs_{\GP} \) yields 
\begin{EQA}
	\fG(\upsvs) - \fG(\upsvs_{\GP})
	&=&
	\frac{1}{2} \| \IFN_{\GP}^{1/2} (\upsvs - \upsvs_{\GP}) \|^{2}
	=
	\frac{1}{2} \| \IFN_{\GP}^{-1/2} \GP^{2} \upsvs \|^{2} 
\label{8chuctc44wckvcuedjequ}
\end{EQA}
and the assertion follows.
\end{proof}

Now we turn to the general case with \( \fs \) satisfying \nameref{LLsT3ref}.
Define
\begin{EQA}
	\IFN_{\GP}
	\eqdef
	\nabla^{2} \fG(\upsvs),
	\qquad
	\bias_{\GP}
	& \eqdef & 
	\| \DFN \, \IFN_{\GP}^{-1} \GP^{2} \upsvs \| 
	\, .
\label{fd9dfhy4ye6fuydfrerf}
\end{EQA}

\begin{theorem}
\label{Pbiasgeneric} 
Let 
\begin{mylist}
\item
the minimizer \( \upsvs = \arginf_{\upsv} \fs(\upsv) \) exists and unique;

\item
\( \fG(\upsv) \) be strongly convex and \( \upsvs_{\GP} = \arginf_{\upsv} \fG(\upsv) \);

\item
\( \fs(\upsv) \) satisfy \nameref{LLsT3ref} at \( \upsvs \) with 
some \( \DFN^{2} \), \( \dltwu_{3} \), and \( \rrn = \rr(\upsvs) \);

\item
for \( \bias_{\GP} = \| \DFN \, \IFN_{\GP}^{-1} \GP^{2} \upsvs \| \) and \( \IFN_{\GP} = \nabla^{2} \fs(\upsvs) + \GP^{2} \);
see \eqref{fd9dfhy4ye6fuydfrerf}:
\begin{EQA}[c]
	\DFN^{2} \leq \dmax^{2} \hspm \IFN_{\GP} 
	\, ,
	\qquad
	\rrn \geq \mfrac{4}{3} \bias_{\GP} 
	\, ,
	\qquad
	\dmax^{2} \hspm \dltwu_{3} \, \bias_{\GP} \leq \mfrac{3}{8} 
	\, .
	\qquad
\label{7fjgjgvuvy44erd52f}
\end{EQA}
\end{mylist}
Then 
\begin{EQA}[rcl]
	\| \DFN (\upsvs_{\GP} - \upsvs) \|
	& \leq &
	\mfrac{4}{3} \, \bias_{\GP}
	\, ,
	\\
	\bigl\| \DFN^{-1} \IFN_{\GP} (\upsvs_{\GP} - \upsvs + \IFN_{\GP}^{-1} \GP^{2} \upsvs) \bigr\|
	& \leq &
	\mfrac{8}{9} \, \dltwu_{3} \, \bias_{\GP}^{2}
	\, ,
\label{11ma3eaelDebbgen}
	\\
	\Bigl| 2 \fG(\upsvs_{\GP}) - 2 \fG(\upsvs) + \| \IFN_{\GP}^{-1/2} \GP^{2} \upsvs \|^{2} \Bigr|
	& \leq &
	\mfrac{1}{2} \, \dltwu_{3} \, \bias_{\GP}^{3}
	\, .
\label{7ywsjhd7wjhjdbiui84kje}
\end{EQA}
\iffourG{
Let, in addition, \( \fs(\upsv) \) satisfy \nameref{LLsT4ref} at \( \upsvs \) with 
the same \( \DFN \), \( \rrn \), and some \( \dltwu_{4} \).
Define
\begin{EQA}
	\bvn_{\GP}
	&=&
	- \IFN_{\GP}^{-1} \{ \GP^{2} \upsvs + \nabla \Tens(\IFN_{\GP}^{-1} \GP^{2} \upsvs) \} \, 
\label{8vfjvr43223efryfuweefG}
\end{EQA}
with \( \Tens(\uv) = \mfrac{1}{6} \langle \nabla^{3} \fs(\upsvs), \uv^{\otimes 3} \rangle \) and
\( \nabla \Tens = \mfrac{1}{2} \langle \nabla^{3} \fs(\upsvs), \uv^{\otimes 2} \rangle \).
Then 
\begin{EQA}[c]
	\| \DFN (\bvn_{\GP} + \IFN_{\GP}^{-1} \GP^{2} \upsvs) \|
	\leq 
	\mfrac{1}{2} \dmax^{2} \dltwu_{3} \, \bias_{\GP}^{2} 
	\, ,
\label{iuvchycvf6e64rygh322b}
\end{EQA}
and with \( \rho_{\GP} \eqdef \dmax^{2} \hspm \dltwu_{3} \, \bias_{\GP} \leq \mfrac{3}{8} \) and
\( \rellG \eqdef 2/(1+\sqrt{1-2\rho_{\GP}}) \in [1,4/3] \), it holds
\begin{EQA}[rcl]
    \| \DFN^{-1} \IFN_{\GP} (\upsvs_{\GP} - \upsvs - \bvn_{\GP}) \|
    & \leq &
    \Bigl\{ \mfrac{1}{6} \dltwu_{4} + \mfrac{\rellG^{2}(\rellG+1)}{4} \, \dmax^{2} \dltwu_{3}^{2} \Bigr\} \, \bias_{\GP}^{3} 
    \, .
\label{DGttGtsGDGm13rGa4b}
\end{EQA}
With \( C(\rellG) \eqdef \mfrac{1}{24} \, \rellG^{2}(2 + 2\rellG - \rellG^{2}) \leq \mfrac{52}{243} \)
for \( \rellG \in [1,4/3] \), 
the remainder
\begin{EQA}[c]
	\Reminder_{\GP}
	\eqdef
	\fG(\upsvs_{\GP}) - \fG(\upsvs) 
	+ \frac{1}{2} \| \IFN_{\GP}^{-1/2} \GP^{2} \upsvs \|^{2} 
	+ \Tens(\IFN_{\GP}^{-1} \GP^{2} \upsvs)
	\, 
\label{ihdviu96ukfr34tyfjf}
\end{EQA}
satisfies
\begin{EQA}[rcccl]
    - \Bigl( \mfrac{1}{24} \dltwu_{4} + C(\rellG) \dmax^{2} \dltwu_{3}^{2} \Bigr) \bias_{\GP}^{4} 
    & \leq &
    \Reminder_{\GP} 
    & \leq &
    \mfrac{1}{24} \dltwu_{4} \, \bias_{\GP}^{4} 
    \, .
    \qquad
\label{3d3Af12DGttGa4b}
\end{EQA}
}{}
\end{theorem}

\begin{proof}
Define \( \fn_{\GP}(\upsv) \) by
\begin{EQA}
	\fn_{\GP}(\upsv) 
	& = &
	\fG(\upsv) - \langle \GP^{2} \upsvs, \upsv \rangle .
\label{f7nhw3dghydfte5w35}
\end{EQA}
Subtracting a linear function preserves strong convexity, so \( \fn_{\GP} \) from \eqref{f7nhw3dghydfte5w35} is strongly convex.
Moreover, \( \nabla \fG(\upsvs) = \GP^{2} \upsvs \) yields \( \nabla \fn_{\GP}(\upsvs) = - \GP^{2} \upsvs + \GP^{2} \upsvs = 0 \).
Hence, \( \upsvs = \arginf \fn_{\GP}(\upsv) \) and \( \fG(\upsv) \) is a linear perturbation \eqref{4hbh8njoelvt6jwgf09} 
of \( \fn_{\GP} \) with \( \Av = \GP^{2} \upsvs \).
Now the results follow from Theorem~\ref{PFiWigeneric2} applied 
with \( \fs(\upsv) = \fG(\upsv) - \langle \Av,\upsv \rangle \),
\( \fn(\upsv) = \fG(\upsv) \), \( \Av = \GP^{2} \upsvs \), and \( \nabla^{2} \fs(\upsvs) = \IFN_{\GP} \).
\iffourG{
The fourth-order expansions follow from Theorem~\ref{Pconcgeneric4} and because \( \avn \) from \eqref{8vfjvr43223efryfuweef}
satisfies \( \avn = \bvn_{\GP} \) and 
\( \nabla^{3} \fs(\upsvs) = \nabla^{3} \fG(\upsvs) = \nabla^{3} \fn_{\GP}(\upsvs) \).
}{}
\end{proof}


\Subsection{A smooth additive penalty}
\label{Slinsmooth}
The case of a general smooth penalty \( \pent_{\GP}(\upsv) \) can be studied similarly to the quadratic case.
Denote 
\( \fG(\upsv) \eqdef \fs(\upsv) + \pent_{\GP}(\upsv) \),
\begin{EQA}
	\upsvs 
	&=& 
	\arginf_{\upsv} \fs(\upsv),
	\quad
	\upsvs_{\GP} = \arginf_{\upsv} \fG(\upsv) 
	=
	\arginf_{\upsv} \bigl\{ \fs(\upsv) + \pent_{\GP}(\upsv) \bigr\}.
\label{8cvkfc9fujf6jnmcer4pen}
\end{EQA}
We study the bias \( \upsvs_{\GP} - \upsvs \) induced by this penalization.
The statement of Theorem~\ref{Pbiasgeneric} and its proof can be easily extended to this situation.
Define
\begin{EQA}
	\IFN_{\GP}
	\eqdef
	\nabla^{2} \fG(\upsvs)
	\, ,
	\quad
	\AvmGP \eqdef \nabla \pent_{\GP}(\upsvs) 
	\, ,
	\quad
	\bias_{\GP}
	& \eqdef & 
	\| \DFN \, \IFN_{\GP}^{-1} \AvmGP \| 
	\, .
\label{fd9dfhy4ye6fuydfrerfG}
\end{EQA} 


\begin{theorem}
\label{Pbiaspen} 
Let  
\begin{mylist}
\item
the minimizer \( \upsvs = \arginf_{\upsv} \fs(\upsv) \) exist and unique;

\item
\( \fG(\upsv) = \fs(\upsv) + \pent_{\GP}(\upsv) \) be strongly convex and \( \upsvs_{\GP} = \arginf_{\upsv} \fG(\upsv) \);

\item
\( \fG(\upsv) \) satisfy \nameref{LLsT3ref} at \( \upsvs \) with 
some \( \DFN^{2} \), \( \dltwu_{3} \), and \( \rrn = \rr(\upsvs) \);

\item
for 
\( \IFN_{\GP} = \nabla^{2} \fG(\upsvs) \),
\( \AvmGP = \nabla \pent_{\GP}(\upsvs) \), and
\( \bias_{\GP} = \| \DFN \, \IFN_{\GP}^{-1} \AvmGP \| \) from \eqref{fd9dfhy4ye6fuydfrerfG}, 
\begin{EQA}[c]
	\DFN^{2} \leq \dmax^{2} \hspm \IFN_{\GP} \, ,
	\qquad
	\rrn
	\geq 
	\mfrac{4}{3} \bias_{\GP} \, ,
	\qquad
	\dmax^{2} \hspm \dltwu_{3} \, \bias_{\GP} \leq \mfrac{3}{8} 
	\, .
	\qquad
\label{8difiyfc54wrbosT4GG}
\end{EQA}
\end{mylist}
Then 
\begin{EQA}[rcl]
\label{odf6fdyr6e4deuewjugG}
	\| \DFN (\upsvs_{\GP} - \upsvs) \|
	& \leq &
	\mfrac{4}{3} \bias_{\GP}
	\, ,
	\\
	\bigl\| \DFN^{-1} \IFN_{\GP} (\upsvs_{\GP} - \upsvs + \IFN_{\GP}^{-1} \AvmGP ) \bigr\|
	& \leq &
	\mfrac{8}{9} \, \dltwu_{3} \, \bias_{\GP}^{2}
	\, ,
	\\
	\Bigl| 2 \fG(\upsvs_{\GP}) - 2 \fG(\upsvs) + \| \IFN_{\GP}^{-1/2} \AvmGP  \|^{2} \Bigr|
	& \leq &
	\mfrac{1}{2} \, \dltwu_{3} \, \bias_{\GP}^{3}
	\, .
\label{7ywsjhd7wjhjdbiui84pen}
\end{EQA}
\iffourG{
Let, in addition, \( \fG(\upsv) \) satisfy \nameref{LLsT4ref} at \( \upsvs \) with 
the same \( \DFN \), \( \rrn \), and some \( \dltwu_{4} \).
With \( \Tens_{\GP}(\uv) = \frac{1}{6} \langle \nabla^{3} \fG(\upsvs), \uv^{\otimes 3} \rangle \),
\( \nabla \Tens_{\GP} = \frac{1}{2} \langle \nabla^{3} \fG(\upsvs), \uv^{\otimes 2} \rangle \), define
\begin{EQA}
	\bvn_{\GP}
	&=&
	- \IFN_{\GP}^{-1} \{ \AvmGP  + \nabla \Tens_{\GP}(\IFN_{\GP}^{-1} \AvmGP ) \} 
	\, .
\label{8vfjvr43223efryfuweefG}
\end{EQA}
Then with \( \rho_{\GP} \eqdef \dmax^{2} \hspm \dltwu_{3} \, \bias_{\GP} \in [0,3/8] \),
and  \( \rellG \eqdef 2/(1+\sqrt{1-2\rho_{\GP}}) \in [1,4/3] \), it holds 
\begin{EQA}
    \| \DFN^{-1} \IFN_{\GP} (\upsvs_{\GP} - \upsvs - \bvn_{\GP}) \|
    & \leq &
    \Bigl\{ \mfrac{1}{6} \dltwu_{4} + \mfrac{\rellG^{2}(\rellG+1)}{4} \, \dmax^{2} \dltwu_{3}^{2} \Bigr\} \, \bias_{\GP}^{3} 
    \, .
    \qquad
\label{DGttGtsGDGm13rGa4bpen}
\end{EQA}
The value remainder \( \Reminder_{\GP} \) given by
\begin{EQA}[c]
	\Reminder_{\GP}
	\eqdef
	\fG(\upsvs_{\GP}) - \fG(\upsvs) 
	+ \frac{1}{2} \| \IFN_{\GP}^{-1/2} \AvmGP \|^{2} 
	+ \Tens_{\GP}(\IFN_{\GP}^{-1} \AvmGP)
	\, ,
\label{ihdviu96ukfr34tyfjf}
\end{EQA}
satisfies with \( C(\rell) \) from \eqref{3d3Af12DGttGa4T} 
\begin{EQA}[rcccl]
    - \Bigl( \mfrac{1}{24} \, \dltwu_{4} + C(\rellG) \dmax^{2} \dltwu_{3}^{2} \Bigr) \bias_{\GP}^{4} 
    & \leq &
    \Reminder_{\GP}
    & \leq &
    \mfrac{1}{24} \dltwu_{4} \, \bias_{\GP}^{4} 
    \, .
\label{3d3Af12DGttGa4bpen}
\end{EQA}
}{}
\end{theorem}

\begin{proof}
Consider \( \fn_{\GP}(\upsv) = \fG(\upsv) - \langle \nabla \pent_{\GP}(\upsvs), \upsv \rangle \).
Then \( \fn_{\GP} \) is strongly convex and \( \nabla \fn_{\GP}(\upsvs) = 0 \) yielding
\( \upsvs = \arginf_{\upsv} \fn_{\GP}(\upsv) \).
Also, \( \fG(\upsv) \) is a linear perturbation of \( \fn_{\GP}(\upsv) \) with 
\( \Av = \AvmGP  = \nabla \pent_{\GP}(\upsvs) \).
Now we can apply Theorem~\ref{PFiWigeneric2}, and under the additional fourth-order assumption Theorem~\ref{Pconcgeneric4}, 
with base objective \( \fn_{\GP} \), perturbed objective \( \fG \), perturbation vector \( \Av = \AvmGP \), and base Hessian \( \IFN_{\GP} \).
\end{proof}


\def\biaslp{\biasl_{\circ}}
\def\Reminderp{\Reminder_{\circ}}
\def\rellp{\rell_{\circ}}
\def\rhop{\rho_{\circ}}
\def\fp{\fs_{\circ}}
\def\dmaxp{\dmax_{\circ}}
\def\upsvp{\upsv_{\circ}}
\def\Tensp{\Tens_{\circ}}
\def\deltp{\Delta}

\Subsection{A generic smooth perturbation}
\label{Sgensmooth}
This section presents a generic perturbation lemma.
Let \( \fp(\upsv) \) be a perturbed version of \( \fs(\upsv) \).
We consider in parallel two optimization problems 
\begin{EQA}
	\upsvs 
	&=& 
	\arginf_{\upsv} \fs(\upsv),
	\qquad
	\upsvp = \arginf_{\upsv} \fp(\upsv) 
	\, .
\label{8cvkfc9fujf6jnmcer4gen}
\end{EQA}
The goal is to describe the displacement \( \upsvp - \upsvs \) induced by this perturbation.
The statement of Theorem~\ref{PFiWigeneric2} and its proof can be easily extended to this situation.
Define
\begin{EQ}[rclrcl]
	\IFN
	& \eqdef &
	\nabla^{2} \fs(\upsvs)
	\, ,
	\qquad
	&
	\IFP
	& \eqdef &
	\nabla^{2} \fp(\upsvs),
	\\
	\biasl
	& \eqdef & 
	\| \DFN \, \IFN^{-1} \nabla \fp(\upsvs) \| 
	\, ,
	\qquad
	&
	\biaslp
	& \eqdef & 
	\| \DFN \, \IFP^{-1} \nabla \fp(\upsvs) \| 
	\, .
\label{fd9dfhy4ye6fuydfrerfGgen}
\end{EQ} 


\begin{theorem}
\label{Pbiasgen} 
Let 
\begin{mylist}
\item
the minimizer \( \upsvs = \arginf_{\upsv} \fs(\upsv) \) exist and unique;

\item
\( \fp(\upsv) \) be strongly convex and \( \upsvp = \arginf_{\upsv} \fp(\upsv) \);

\item
\( \fp(\upsv) \) satisfy \nameref{LLsT3ref} at \( \upsvs \) with 
some \( \DFN^{2} \), \( \dltwu_{3} \), and \( \rrn = \rr(\upsvs) \);

\item
for 
\( \IFP = \nabla^{2} \fp(\upsvs) \),
\( \biaslp = \| \DFN \, \IFP^{-1} \nabla \fp(\upsvs) \| \) from \eqref{fd9dfhy4ye6fuydfrerfGgen}, and some \( \dmaxp > 0 \)
\begin{EQA}[c]
	\DFN^{2} \leq \dmaxp^{2} \hspm \IFP \, ,
	\qquad
	\rrn
	\geq 
	\mfrac{4}{3} \biaslp \, ,
	\qquad
	\dmaxp^{2} \hspm \dltwu_{3} \, \biaslp \leq \mfrac{3}{8} 
	\, .
	\qquad
\label{8difiyfc54wrbosT4Gg}
\end{EQA}
\end{mylist}
Then 
\begin{EQ}[rcl]
\label{odf6fdyr6e4deuewjugGgen}
	\| \DFN (\upsvp - \upsvs) \|
	& \leq &
	\mfrac{4}{3} \biaslp
	\, ,
	\\
	\bigl\| \DFN^{-1} \IFP (\upsvp - \upsvs + \IFP^{-1} \nabla \fp(\upsvs) ) \bigr\|
	& \leq &
	\mfrac{8}{9} \, \dltwu_{3} \, \biaslp^{2}
	\, ,
	\\
	\Bigl| 2 \fp(\upsvp) - 2 \fp(\upsvs) + \| \IFP^{-1/2} \nabla \fp(\upsvs) \|^{2} \Bigr|
	& \leq &
	\mfrac{1}{2} \, \dltwu_{3} \, \biaslp^{3}
	\, .
\end{EQ}
\iffourG{
If, also, \( \fp(\upsv) \) satisfies \nameref{LLsT4ref}  at \( \upsvs \) with 
the same \( \DFN \), \( \rrn \), and some \( \dltwu_{4} \)
then with \( \Tensp(\uv) = \frac{1}{6} \langle \nabla^{3} \fp(\upsvs), \uv^{\otimes 3} \rangle \),
\( \nabla \Tensp = \frac{1}{2} \langle \nabla^{3} \fp(\upsvs), \uv^{\otimes 2} \rangle \), and
\begin{EQA}
	\bvp
	&=&
	- \IFP^{-1} \{ \nabla \fp(\upsvs)  + \nabla \Tensp(\IFP^{-1} \nabla \fp(\upsvs) ) \} \, ,
\label{8vfjvr43223efryfuweefG}
\end{EQA}
it holds 
\begin{EQA}[rcl]
    \| \DFN^{-1} \IFP (\upsvp - \upsvs - \bvp) \|
    & \leq &
    \Bigl\{ \mfrac{1}{6} \dltwu_{4} + \mfrac{\rellp^{2}(\rellp+1)}{4} \, \dmaxp^{2} \dltwu_{3}^{2} \Bigr\} \biaslp^{3} 
    \, ,
\label{3d3Af12DGttGa4bgen}
\end{EQA}
where \( \rellp \eqdef 2/(1+\sqrt{1-2\rhop}) \in [1,\mfrac{4}{3}] \)  
for \( \rhop \eqdef \dmaxp^{2} \hspm \dltwu_{3} \, \biaslp \in [0,\mfrac{3}{8}] \).
The value remainder 
\begin{EQA}[c]
	\Reminderp
	\eqdef
	\fp(\upsvp) - \fp(\upsvs) 
    + \frac{1}{2} \| \IFP^{-1/2} \nabla \fp(\upsvs) \|^{2} + \Tensp(\IFP^{-1} \nabla \fp(\upsvs) )
    \, 
\end{EQA}
fulfills with \( C(\rell) \) from \eqref{3d3Af12DGttGa4T}
\begin{EQA}[rcl]
    - \Bigl( \mfrac{1}{24} \, \dltwu_{4} + C(\rellp) \dmaxp^{2} \dltwu_{3}^{2} \Bigr) \biaslp^{4} 
    & \leq &
    \Reminderp
    \leq 
    \mfrac{1}{24} \, \dltwu_{4} \, \biaslp^{4} 
    \, .
    \qquad
\label{3d3Af12DGttGa4bgen}
\end{EQA}
}{}
\end{theorem}

\begin{proof}
Consider \( \fh(\upsv) = \fp(\upsv) - \langle \nabla \fp(\upsvs), \upsv \rangle \).
Then \( \fh \) is strongly convex and \( \nabla \fh(\upsvs) = 0 \) yielding
\( \upsvs = \arginf_{\upsv} \fh(\upsv) \).
Also, \( \fp(\upsv) \) is a linear perturbation of \( \fh(\upsv) \) with 
\( \Av = \AvmP  = \nabla \fp(\upsvs) \); see \eqref{4hbh8njoelvt6jwgf09}.
Now all the statements of Theorem~\ref{PFiWigeneric2} 
\iffourG{and Theorem~\ref{Pconcgeneric4}}{}
apply to \( \upsvp \).
\end{proof}

In some settings, using the Hessian \( \IFP \) is inconvenient.
The next result restates the third-order bounds with \( \IFN \) in place of \( \IFP \)
under a condition on \( \IFP - \IFN \).

\begin{theorem}
\label{TbiasgenP}
Let
\begin{mylist}
\item
\( \fs(\upsv) \) be strongly convex and \( \upsvs = \arginf_{\upsv} \fs(\upsv) \),
\( \IFN = \nabla^{2} \fs(\upsvs) \);

\item
\( \fp(\upsv) \) be strongly convex and \( \upsvp = \arginf_{\upsv} \fp(\upsv) \),
\( \IFP = \nabla^{2} \fp(\upsvs) \);

\item
\( \fp(\upsv) \) satisfy \nameref{LLsT3ref} at \( \upsvs \) with 
some \( \DFN^{2} \), \( \dltwu_{3} \), and \( \rrn = \rr(\upsvs) \);

\item 
For some \( \deltp \)
\begin{EQA}[c]
	\bigl\| \DFN^{-1} ( \IFP - \IFN ) \DFN^{-1} \bigr\|
	\leq 
	\deltp
	\, .
\label{udhiiiyt445ghgkfgrr}
\end{EQA}

\item
for 
\( \biasl = \| \DFN \, \IFN^{-1} \nabla \fp(\upsvs) \| \) from \eqref{fd9dfhy4ye6fuydfrerfGgen}, and some \( \dmax > 0 \)
with \( \dmax^{2} \deltp < 1 \)
\begin{EQA}[c]
	\DFN^{2} \leq \dmax^{2} \hspm \IFN 
	\, ,
	\quad
	\rrn
	\geq 
	\mfrac{4}{3 (1 - \dmax^{2} \deltp)} \, \biasl
	\, ,
	\quad
	\mfrac{\dmax^{2}}{(1 - \dmax^{2} \deltp)^{2}} \hspm \dltwu_{3} \, \biasl \leq \mfrac{3}{8} 
	\, .
	\qquad
\label{8difiyfc54wrbosT4Ggg}
\end{EQA}
\end{mylist}
Then with \( \biasls \eqdef \mfrac{\biasl}{1 - \dmax^{2} \deltp} \)
\begin{EQ}[rcl]
\label{odf6fdyr6e4deuewjugGgg}
	\| \DFN (\upsvp - \upsvs) \|
	& \leq &
	\mfrac{4}{3} \, \biasls
	\, ,
	\\
	\bigl\| \DFN^{-1} \IFN (\upsvp - \upsvs + \IFN^{-1} \nabla \fp(\upsvs) ) \bigr\|
	& \leq &
	\mfrac{8}{9} \, \dltwu_{3} \, \biasls^{2}
	+ \mfrac{4}{3} \, \deltp \, \biasls
	\, ,
	\qquad
	\qquad
	\\
	\Bigl| 2 \fp(\upsvp) - 2 \fp(\upsvs) + \| \IFN^{-1/2} \nabla \fp(\upsvs) \|^{2} \Bigr|
	& \leq &
	\mfrac{1}{2} \, \dltwu_{3} \, \biasls^{3}
	+ \deltp \, \biasls^{2}
	\, .
	\qquad
\end{EQ}
\end{theorem}

\begin{proof}
Condition \eqref{udhiiiyt445ghgkfgrr} implies \( \IFP \geq \IFN - \deltp \DFN^{2} \)
and it follows from \( \DFN^{2} \leq \dmax^{2} \IFN \) that
\begin{EQA}[rcl]
	\DFN^{2} 
	& \leq &
	\dmaxp^{2} \IFP 
	\, ,
	\qquad
	\dmaxp^{2} \eqdef \mfrac{\dmax^{2}}{1 - \dmax^{2} \deltp} 
	\, .
\end{EQA}
Also, for any \( \uv \in \R^{\dimp} \)
\begin{EQA}[rcl]
	&& \nquad
	\| \DFN (\IFP^{-1} - \IFN^{-1}) \uv \|
	=
	\| \DFN \IFP^{-1} (\IFP - \IFN) \IFN^{-1} \uv \|
	\\
	&=&
	\| \DFN \IFP^{-1} \DFN \; \DFN^{-1} (\IFP - \IFN) \DFN^{-1} \; \DFN \IFN^{-1} \uv \|
	\\
	& \leq &
	\| \DFN \IFP^{-1} \DFN \| \; \| \DFN^{-1} (\IFP - \IFN) \DFN^{-1} \| \;  
	\| \DFN \IFN^{-1} \uv \|
	\leq 
	\mfrac{\dmax^{2}\deltp}{1 - \dmax^{2} \deltp} \, \| \DFN \IFN^{-1} \uv \|
	\, 
\end{EQA}
yielding
\begin{EQA}[c]
	\| \DFN \IFP^{-1} \uv \|
	\leq 
	\mfrac{1}{1 - \dmax^{2} \deltp} \, \| \DFN \IFN^{-1} \uv \|
	\, .
\label{jyuiujere334rtfyee}
\end{EQA}
For \( \uv = \nabla \fp(\upsvs) \), this implies with
\( \biasl = \| \DFN \, \IFN^{-1} \nabla \fp(\upsvs) \| \),
\( \biaslp = \| \DFN \, \IFP^{-1} \nabla \fp(\upsvs) \| \) that 
\begin{EQA}[rcl]
	|\biaslp - \biasl|
	& \leq &
	\| \DFN (\IFP^{-1} - \IFN^{-1}) \nabla \fp(\upsvs) \|
	\leq 
	\mfrac{\dmax^{2} \deltp}{1 - \dmax^{2} \deltp} \, \biasl
	\, ,
\end{EQA}
and thus
\(	\biaslp \leq \mfrac{1}{1 - \dmax^{2} \deltp} \, \biasl 
	=
	\biasls
\).
Hence, \eqref{8difiyfc54wrbosT4Ggg} ensures \eqref{8difiyfc54wrbosT4Gg}.
This justifies the claims in \eqref{odf6fdyr6e4deuewjugGgen}.
In particular,
\( \| \DFN (\upsvp - \upsvs) \| \leq \mfrac{4}{3} \, \deltp \, \biaslp \leq \mfrac{4}{3} \, \deltp \, \biasls \).
Further,
\begin{EQA}[rcl]
	&& \nquad
	\bigl\| \DFN^{-1} \IFP (\upsvp - \upsvs + \IFP^{-1} \nabla \fp(\upsvs) )
	- \DFN^{-1} \IFN (\upsvp - \upsvs + \IFN^{-1} \nabla \fp(\upsvs) ) \bigr\|
	\\
	&=&
	\bigl\| \DFN^{-1} (\IFP - \IFN) (\upsvp - \upsvs) \bigr\|
	\leq 
	\| \DFN^{-1} (\IFP - \IFN) \DFN^{-1} \| \; \| \DFN (\upsvp - \upsvs) \|
	\\
	& \leq &
	\mfrac{4}{3} \, \deltp \, \biaslp
	\leq 
	\mfrac{4}{3} \, \deltp \, \biasls
	\, .
\end{EQA}
This and \eqref{odf6fdyr6e4deuewjugGgen} imply the second bound in \eqref{odf6fdyr6e4deuewjugGgg}.
Also, for \( \uv \in \R^{\dimp} \), by \eqref{udhiiiyt445ghgkfgrr} and \eqref{jyuiujere334rtfyee}
\begin{EQA}[rcl]
	&& \nquad
	\bigl| \| \IFP^{-1/2} \uv \|^{2} - \| \IFN^{-1/2} \uv \|^{2} \bigr|
	=
	\bigl| \uv^{\T} (\IFP^{-1} - \IFN^{-1}) \uv \bigr|
	=
	\bigl| \uv^{\T} \IFN^{-1} (\IFP - \IFN) \IFP^{-1} \uv \bigr|
	\\
	&=&
	\bigl| (\DFN \IFN^{-1} \uv)^{\T} \DFN^{-1} (\IFP - \IFN) \DFN^{-1} (\DFN \IFP^{-1} \uv) \bigr|
	\leq 
	\deltp \| \DFN \IFN^{-1} \uv \| \; \| \DFN \IFP^{-1} \uv \|
	\\
	& \leq &
	\mfrac{1}{1 - \dmax^{2} \deltp} \, \deltp \| \DFN \IFN^{-1} \uv \|^{2}
	\, .
\end{EQA}
With \( \uv = \nabla \fp(\upsvs) \), this and \eqref{odf6fdyr6e4deuewjugGgen}
imply the final bound in \eqref{odf6fdyr6e4deuewjugGgg}.
\end{proof}

\appendix
\Section{Technical assertions}
This section collects a few lemmas used in the main proofs.

\begin{lemma}
\label{Llin23}
Let \( \AFN \geq \dmax^{-2} \Id_{\dimp} \) and let \( \afv \in \R^{\dimp} \) satisfy
\( \rho \eqdef \dmax^{2} \dltwu \| \afv \| \leq \mfrac{1}{2} \) for \( \dltwu \geq 0 \).
Define \( \rell \eqdef \mfrac{2}{1 + \sqrt{1 - 2 \rho}} \) and let
\( \rr \geq \rell \| \afv \| \) but \( \dmax^{2} \dltwu \, \rr \leq 1 \).
Then
\begin{EQA}
\label{hdcf6wyheuv76e34r35eycv}
	\sup_{\| \uv \| \leq \rr} \Bigl( \mfrac{\dltwu}{3} \| \uv \|^{3} - (\uv - \afv)^{\T} \AFN (\uv - \afv) \Bigr)
	& \leq & 
	\mfrac{\dltwu \rell^{2}(3 - \rell)}{6} \, \| \afv \|^{3} \, ,
	\\
	\inf_{\| \uv \| \leq \rr} \Bigl( \mfrac{\dltwu}{3} \| \uv \|^{3} + (\uv - \afv)^{\T} \AFN (\uv - \afv) \Bigr)
	& \leq & 
	\mfrac{\dltwu}{3} \, \| \afv \|^{3} \, .
\label{hdcf6wyheuv76e34r35eycvm}
\end{EQA}
\end{lemma}

\begin{proof}
For \( \afv = 0 \), the result is obvious.
As \( \AFN \geq \dmax^{-2} \Id_{\dimp} \), it holds with \( b \eqdef \| \afv \| > 0 \)
\begin{EQA}
	&& \nquad
	\sup_{\| \uv \| \leq \rr} \Bigl( \mfrac{\dltwu}{3} \| \uv \|^{3} - (\uv - \afv)^{\T} \AFN (\uv - \afv) \Bigr)
	\leq 
	\sup_{\| \uv \| \leq \rr} \Bigl( \mfrac{\dltwu}{3} \| \uv \|^{3} - \dmax^{-2} \| \uv - \afv \|^{2} \Bigr)
	\\
	& \leq & 
	\sup_{\| \uv \| \leq \rr} \; \sup_{\| \xv \| = b} \Bigl( \mfrac{\dltwu}{3} \| \uv \|^{3} - \dmax^{-2} \| \uv - \xv \|^{2} \Bigr)
	\\
	&=&
	\sup_{\| \uv \| \leq \rr} \Bigl( \mfrac{\dltwu}{3} \| \uv \|^{3} - \dmax^{-2} (\| \uv \| - b)^{2} \Bigr)
	=
	b^{2} \dmax^{-2}  \sup_{t \in [0,\rr/b]} \Bigl( \mfrac{\dmax^{2} \dltwu b}{3} t^{3} - (t - 1)^{2} \Bigr)
	\, .
\end{EQA}
By Lemma~\ref{Ltt3tb2}, under \( \dmax^{2} \dltwu b \leq 1/2 \) and \( \rr/b \geq \rell \eqdef \mfrac{2}{1 + \sqrt{1 - 2 \dltwu b}} \),
\begin{EQA}[c]
	b^{2} \dmax^{-2} \sup_{t \in [0,\rr/b]} \Bigl( \mfrac{\dmax^{2} \dltwu b}{3} t^{3} - (t - 1)^{2} \Bigr)
	=
	\mfrac{\dltwu b^{3} \rell^{2}(3 - \rell)}{6} 
	\, .
\end{EQA}
Further, 
the minimum of the left-hand side of \eqref{hdcf6wyheuv76e34r35eycvm} is not larger than the value
at the point \( \uv = \afv \).
This completes the proof.
\end{proof}

\begin{lemma}
\label{Ltt3tb2}
Let \( \rho \in [0,1/2] \).
Define \( \rell \eqdef \mfrac{2}{1 + \sqrt{1 - 2 \rho}} \in [1,2] \). 
Then the function 
\begin{EQA}[rcl]
	\phi(t)
	& \eqdef &
	\mfrac{\rho}{3} t^{3} - (t - 1)^{2}
\end{EQA} 
satisfies for any \( T \in [\rell,1/\rho] \)
\begin{EQA}[rcl]
	\sup_{t \in [0,T]} \phi(t)
	&=&
	\phi(\rell)
	= 
	\mfrac{\rho \rell^{2} (3 - \rell)}{6}
	\, .
\label{khvdfi8yud8w2vouvww}
\end{EQA}
\end{lemma}

\begin{proof}
The case \( \rho = 0 \) is straightforward.
Consider \( \rho > 0 \).
As \( \phi''(t) = 2 \rho t - 2 \leq 0 \) for \( t \leq 1/\rho \),
the function \( \phi \) is concave on \( [0,T] \) for \( T \leq 1/\rho \).
The stationary point of \( \phi \) solves \( \rho t^{2} - 2 (t-1) = 0 \). 
This leads to 
\begin{EQA}[rcl]
	t 
	&=&
	\mfrac{1 - \sqrt{1 - 2 \rho}}{\rho}
	=
	\mfrac{2}{1 + \sqrt{1 - 2 \rho}}
	=
	\rell
	\leq 
	T
	\, .
\end{EQA}
Simple algebra yields \( \rell - 1 = \rell (1 - \sqrt{1 - 2 \rho})/2 = \rho \rell^{2}/2 \) and \( t = \rell \) yields
\begin{EQA}[rcl]
	\mfrac{\rho}{3} t^{3} - (t - 1)^{2}
	& = &
	\mfrac{\rho}{3} \rell^{3} - (\rell - 1)^{2}
	=
	\mfrac{\rho}{3} \rell^{3} - \mfrac{\rho}{2} (\rell - 1) \rell^{2}
	=
	\mfrac{\rho \rell^{2} (3 - \rell)}{6} 
	\, ,
\end{EQA}
and \eqref{khvdfi8yud8w2vouvww} follows.
\end{proof}

\begin{lemma}
\label{L3poly}
Let a function \( P(t) \) be two times continuously differentiable and, given \( T > 0 \), 
satisfy
\( P''(t) \leq 0 \) for \( t \in [0,T] \) and \( P'(T) = 0 \).
Then 
\begin{EQA}[c]
	P(t) \leq P(T)
	\, , 
	\quad
	t \in [0,T]
	\, .
\end{EQA}
\end{lemma}

\begin{proof}
Condition \( P''(t) \leq 0 \) for \( t \in [0,T] \) implies that
\( P'(t) \) is decreasing on \( \II \).
Further, the edge conditions \( P'(T) = 0 \) yields \( P'(t) \geq 0 \) for \( t \in [0,T] \).
Therefore, \( P(t) \) is increasing on \( [0,T] \) and satisfies \( P(t) \leq P(T) \).
\end{proof}

\section*{Conclusion}
The paper systematically studies the effect of changing the objective function by a linear, quadratic, or smooth perturbation.
The obtained results provide careful expansions for the solution and the value of the perturbed optimization problem.
These expansions can be used as building blocks in different areas, including statistics and machine learning, 
quasi-Newton optimization, and uncertainty quantification for inverse problems, among many others. 

\bibliographystyle{apalike}

\bibliography{exp_ts,listpubm-with-url}

\end{document}